\documentclass[%
 reprint,
 amsmath,amssymb,
 aps,
floatfix
]{revtex4-2}

\usepackage{graphicx}
\usepackage{dcolumn}
\usepackage{bm}
\usepackage{hyperref}
\usepackage{amsmath}
\usepackage{enumerate}
\usepackage{stfloats}
\usepackage{enumitem}
\usepackage{subcaption}
\usepackage{xcolor}
\usepackage{float}
\usepackage[skip=0.5ex]{subcaption}
\usepackage{lipsum}

\hypersetup{
    colorlinks=true,
    linkcolor=blue,
     urlcolor=blue,
     citecolor=blue,
    pdftitle={Overleaf Example},
    pdfpagemode=FullScreen,
    }
    
\begin{document}

\preprint{APS/123-QED}

\title{
Schr\"odinger Bridges over Kinetic Swarming Models
}

\author{Asmaa Eldesoukey}\email{asmaae@iastate.edu}
\affiliation{%
 Department of Aerospace Engineering, Iowa State University, Ames, IA 50011, USA}%

\author{Md Zulfiqur Haider}%
 \email{zulfiqur@iastate.edu}
\affiliation{
Department of Mathematics, Iowa State University, Ames, IA 50011, USA
}%

\author{Italo Napolitano}
\email{i.napolitano@ssmeridionale.it}
\affiliation{%
Modeling and Engineering Risk and Complexity Department, Scuola
Superiore Meridionale, via Mezzocannone 4, 80138, Naples, Italy
}%

\author{Yongxin Chen}%
\email{yongchen@gatech.edu}
\affiliation{
Yongxin Chen is with the School of Aerospace Engineering, Georgia Institute of Technology, Atlanta, GA 30332, USA
}%

\author{Abhishek Halder}
\email{ahalder@iastate.edu}
\affiliation{%
Department of Aerospace Engineering, Iowa State University, Ames, IA 50011, USA
}%
\affiliation{%
Department of Applied Mathematics, University of California Santa Cruz, 1156 High Street, Santa Cruz, CA 95064, USA
}%

\date{\today}

\begin{abstract}
Paradigmatic interaction models explain how collective behaviors can emerge in complex systems from interactions among the constituent agents. In bio-inspired swarms, however, interactions alone may not suffice to bring the population to a desired aggregate configuration within a prescribed time horizon, as needed in applications ranging from targeted therapy to collective transport and emergency evacuation. In the present work, we consider finite-horizon minimum-energy collective steering for inertial swarms that are subject to stochastic disturbances. We focus on the mean-field representations of these multi-agent systems driven by Cucker--Smale alignment or Morse attraction--repulsion interactions. Our objective is to steer the swarm between prescribed endpoint distributions using a state-feedback control, where the endpoint specifications can be full phase-space distributions (positions and velocities) or position marginals alone. Our formalism is rooted in the theory of Schr\"odinger bridges, which has inspired contemporary developments spanning statistical inference, biological modeling, stochastic control, and generative learning. Within the bridges framework, the uncontrolled interacting stochastic dynamics are viewed as a prior model, and the optimal control as the minimum-energy corrective drift needed to realize the prescribed distributions. We derive nonlinear, coupled necessary optimality systems with a time-symmetric structure reminiscent of classical Schr\"odinger bridges, and propose nested fixed-point schemes to numerically solve them.
Numerical examples show that the obtained optimal control (corrective drift) can dynamically exploit or counteract the interaction forces, depending on whether the latter are favorable or adversarial to the steering task.


\end{abstract}

\maketitle


\section{Introduction}

\emph{Emergent behaviors} refer to large-scale patterns that result from self-organization within complex systems that often operate without centralized leadership \cite{artime2022origin}. Gregarious swarms exhibiting such behaviors are pervasive in nature: birds flock \cite{ballerini2008empirical}, fish school \cite{lopez2012behavioural}, mammals engage in herding \cite{king2012selfish,kameda2015herd}, and algae form blooms \cite{harvey2011avoidance}, all of which garnered significant research interest. Beyond nature, the notion of decentralized self-organization has inspired the development of active systems such as synthetic self-propelled microswimmers \cite{bechinger2016active,bricard2013emergence,ziepke2022multi} and externally actuated microrobotic swarms \cite{ceron2023programmable}. These systems utilize coordination among the constituent agents to accomplish complicated tasks such as targeted therapy \cite{wang2021ultrasound,xie2025microrobotic} and environmental remediation \cite{shivalkar2021autonomous}.

Numerous studies have focused on modeling collective behaviors in biological or biologically inspired swarms, to characterize the interaction mechanisms that give rise to these emergent phenomena. For example, \emph{the Vicsek model} realizes flocking by combining agents' self-propulsion and local alignment, with each agent moving at a constant speed while matching its direction to the average of its neighbors \cite{czirok1999collective,vicsek2012collective}. \emph{The Cucker--Smale model} describes consensus formation within groups via velocity alignment, mediated by distance-dependent communication weights \cite{cucker2007emergent,cucker2007mathematics}. \emph{Morse-potential} driven models lead to spatial organization with a confluence of attractive and repulsive interactions to describe both ensemble cohesion and collision avoidance among agents \cite{d2006self,carrillo2013new}.

These and other swarming models have become paradigmatic frameworks for formalizing how flocking, milling \cite{Carrillo2014}, and coalescence emerge without external control. However, interactions may not solely ensure that the swarm attains a desired aggregate configuration within a prescribed time horizon. This specification is relevant in applications where timely pattern formation or reconfiguration is of interest, such as targeted therapy \cite{qin2023spatiotemporal,zhong2026autonomous} and self-assembly \cite{kim2014spatially,nodozi2023physics,nodozi2023neural}.  

The present work addresses the problem of \emph{finite-horizon minimum-energy collective steering} of swarms subject to noise toward a desired aggregate configuration.
 The term ``collective steering" means that the control command is broadcast to all particles/agents with the aim of regulating the dynamics of the collective, rather than individual trajectories.  Specifically, given (i) a stochastic swarming kinetic model that governs interacting particles and (ii) initial and target probability distributions, we seek an optimal control input that guides the swarm between the prescribed distributions by a chosen terminal time. We focus on \emph{kinetic models} driven by either the Cucker--Smale or Morse potential interactions. We also consider two scenarios in which the prescribed distributions are either over the full phase space (position and velocity) or over positions only. Our framework is rooted in the theory of \emph{Schr\"odinger bridges}, a maximum-likelihood problem originally posed and solved by Erwin Schr\"odinger.

In the 1930s, Schr\"odinger posed a thought experiment \cite{schrodinger1931,schrodinger1932} in his attempt to understand 
quantum mechanics from a statistical standpoint. He considered \emph{non-interacting} particles that follow known macroscopic prior dynamics (in particular, the Fokker-Planck or Kolmogorov's forward equation). The probability distributions of such particles are observed only at the start and end of the experiment. Assuming that the empirical final distribution differs from that predicted by the prior dynamics, he therefore asked: \emph{among all stochastic evolutions consistent with the observed initial and final marginals, which one is the most likely?} The resulting inference (inverse) problem is now known as the \emph{Schr\"odinger bridge problem}. It seeks the most likely paths that particles traverse between the endpoint empirical distributions. 

The Schr\"odinger bridge problem is formalized to identify the most likely probability law on sample paths that aligns with the observed marginals. This law effectively acts as a ``bridge" between the observations. Notably, among all laws that match the prescribed marginals, the most likely is the one that represents the smallest perturbation of the prior law, as measured by the \emph{relative entropy}. Schr\"odinger's solution to this problem defines a bridge whose marginal densities appear as a product of a forward-diffusing function and a backward-diffusing one (see Section~\ref{sec:background}). Interestingly, this structure, i.e., the factorization of the probability density into a product of two functions, bears a ``striking analogy" \cite{chetrite2021schrodinger} to that obtained from his equation in quantum mechanics.

Later, Dai Pra \cite{daipra1991new} established the equivalence between the Schr\"odinger bridge problem and the minimum-energy control problem of steering non-interacting particles between two endpoint marginals using Girsanov's theorem. In this, the optimal control represents the smallest, in the sense of expected energy, required added drift to the prior dynamics so as to match the endpoint marginals (see also \cite{chen2016optimal}). This development endowed the Schr\"odinger bridge problem, beyond the inference perspective, with yet another interpretation as a stochastic optimal control problem. This \emph{inference--control duality} has thus far enabled the application of the Schr\"odinger bridge framework in physically motivated contexts, including self-assembly, stochastic cooling, oscillator synchronization, systems with particle loss, first-passage processes, the inference of potential landscapes and many others; see, e.g., \cite{nodozi2023neural,nodozi2022schrodinger,chen2015fast,chen2022most,eldesoukey2024schrodinger,movilla2024inferring,chizat2022trajectory,teter2025probabilistic}.
\vspace{5pt}\\
\noindent \textbf{Scope and Significance}\\[5pt]
We consider two second-order stochastic swarming models: the Cucker--Smale model and a Morse-potential-driven model.  For each, we formulate a stochastic minimum-energy control problem with prescribed initial and final distributions, either over the full phase space of positions and velocities or over positions alone. We consider these interacting particle systems in the \emph{mean-field limit}, where the number of agents approaches infinity and the collective density evolves according to a \emph{nonlinear Vlasov--Fokker--Planck integro-partial differential equation (PDE)}. 

Our work thereby extends the classical Schr\"odinger bridge framework to \emph{interacting kinetics} under both complete and partial endpoint observations. We derive the corresponding nonlinear, coupled necessary optimality systems and propose novel nested iterative schemes for computing solutions of these systems. Our numerical examples demonstrate how the obtained control may counteract or exploit interactions, which can either increase or reduce the required control effort compared to the non-interacting case, while consistently having lower cost than an interaction-canceling baseline.

Equivalently, from an inference perspective, the \emph{uncontrolled} Vlasov--Fokker--Planck equation may be viewed as a prior evolution. Among all evolutions compatible with the given distributions, the minimum-energy control problem seeks the one requiring the smallest correction to the prior dynamics, characterizing the most probable evolution. The inference--control connection underlying this formulation for weakly interacting inertial systems was established by Fischer \cite{fischer2014form}. However, that work did not consider the associated endpoint-constrained Schr\"odinger bridge problem studied here. 
When distributions are specified only over positions, our formulation pertains to scenarios where only partial information about the system is available, for example, through aggregate snapshots of the swarm, leaving the corresponding velocity degree of freedom to be inferred. 

The practical motivation behind this work is to provide a framework for steering large-scale interacting particles collectively to fulfill desired functionalities within a specified time horizon, such as pattern formation required in microrobotic swarms \cite{ceron2023programmable} for, e.g., targeted therapy \cite{wang2021ultrasound,xie2025microrobotic}. From the inference perspective, within the context of complex systems, the corrective drift may represent \emph{unmodeled forcing}, e.g., arising in response to external social or environmental drivers, such as interventions by authorities or policymakers in a population \cite{piccoli2015control}, or the presence of predators or obstacles in animal groups. Unmodeled dynamics can also represent internal factors, such as an agent's autonomous intent, a differentiated role \cite{cucker2008flocking,shen2008cucker,yuan2022flocking,ren2025collision}, or self-propulsion \cite{lukeman2010inferring}, or additional underlying interaction mechanisms \cite{cucker2010avoiding,park2010cucker}.

\vspace{5pt}
\noindent \textbf{Relation to Prior Work}\\[5pt]
Existing control strategies for interacting-particle systems typically aim to promote consensus, induce flocking, or stabilize selected macroscopic states. For example, sparse and optimal control policies have been designed for Cucker-Smale dynamics to enforce flocking or consensus \cite{piccoli2015control, bailo2018optimal}. For attractive-repulsive swarms, bounded control strategies have been developed for asymptotic stabilization toward flock and mill configurations \cite{carrillo2022controlling}. Bio-inspired feedback approaches such as shepherding have also been proposed to regulate long-time collective configurations \cite{maffettone2022continuification,lama2025nonreciprocal,maffettone2026bio,napolitano2026hierarchical,catello2026sparse}. These approaches generally enforce the desired target behavior through either soft penalties or stabilization objectives. The present work, however, imposes the terminal probability distribution as a hard constraint at a prescribed finite time and seeks to minimize the control effort required to reach it.

In the context of inference, recent studies \cite{heins2024collective,king2025inferring,miller2023learning} have focused on learning interaction mechanisms from time-resolved data such as particle trajectories. In contrast, we infer a collective evolution that reconciles a prior model with aggregate endpoint data. We identify a data-informed corrective drift rather than learning, e.g., the communication weights in the Cucker--Smale interaction kernel.

Previous work \cite{chen2019multi,chiarini2022entropic} integrated kinetic models into the Schr\"odinger bridge framework with spatial endpoint constraints, but without interactions. Backhoff et al.~\cite{backhoff2020mean} extended the classical Schr\"odinger bridge setting to interacting mean-field systems but with first-order dynamics. Orland \cite{orland2025schrodinger} developed a quantum-inspired formulation for interacting inertialess Schr\"odinger bridges, yielding coupled forward and backward single-particle wave functions.
On the numerical side, \cite{chen2024density,eldesoukey2026generalized} developed methods for solving such first-order mean-field bridge problems. 
Other related work has also considered minimum-energy steering for finite populations of second-order stochastic oscillators \cite{nodozi2023neural}, and has developed controlled nonlocal mean-field models for interacting particle populations in collective assembly \cite{nodozi2023controlled}. However, the former does not consider the kinetic mean-field limit, and the latter does not formulate an endpoint-constrained Schr\"odinger bridge problem.

Our work bridges the gap between kinetic and interacting Schr\"odinger bridge formulations, addressing kinetic systems with nonlocal mean-field interactions and formulating finite-horizon Schr\"odinger bridge problems for both prescribed phase-space endpoint distributions and spatial marginals. We develop the corresponding optimality structure for Cucker--Smale and Morse-type interactions, and propose fixed-point schemes for their computation.

\section{Schr\"odinger Bridge Formulations}\label{sec:background}
\subsection{Classical Formulation}\label{sec:classical-LD}
Consider a collection of particles (agents) in motion over the time window $[0,T]$. Each particle is governed by the It\^o diffusion
\begin{align}\label{eq:Ito}
    {\rm d}X_t = b_t(X_t) \,{\rm d}t + \sigma \, {\rm d}B_t, \quad X_0 \sim q_0(x),
\end{align}
where $X_t \in\mathbb R^d$.
The noise intensity is denoted by $\sigma>0$, which we take to be constant for simplicity, and $B_t$ is a $ d$-dimensional standard Brownian motion. Initially, particles are distributed according to a known probability density $q_0 = \rho_{\rm initial}$. Their density $q_t$ evolves according to the Fokker--Planck equation
\begin{align}\label{eq:FP}
    \partial_t q_t + \nabla_x \cdot (b_t q_t) = \frac{\sigma^2}{2} \Delta_x q_t. 
\end{align}
Equivalently, if $k(s,x,t,y)$ denotes the Markov transition kernel (propagator) from $x$ at instant $s$ to $y$ at instant $t$ that corresponds to \eqref{eq:FP}, we write
\begin{align*}
    q_t(y) = \int_{\mathbb R^d} k(0,x,t,y) q_0(x)  \, {\rm d}x,
\end{align*}
for any $y \in \mathbb R^d$. 

The sample paths of the It\^o process in \eqref{eq:Ito} are continuous functions on the interval $[0, T]$.
We denote by $Q$ the probability law of these paths, assigning probabilities to sets of trajectories.  
Suppose, however, that at the final time, the particles are observed to have a probability density $\rho_{\rm final}$, where 
\begin{align}
   \rho_{\rm final}(y)  \neq  \int_{\mathbb R^d} k(0,x,T,y) q_0(x)  \, {\rm d}x.
\end{align}
Such an observation is considered \emph{atypical} if it differs significantly from what is predicted by the 
prior model 
described by \eqref{eq:Ito} or \eqref{eq:FP}, and the corresponding law $Q$.

Among all candidate path laws consistent with the atypical observation, the Schr\"odinger bridge problem seeks the most likely alternative to the prior $Q$. 
A seminal result in the large deviations theory, Sanov's theorem \cite{sanov1958probability}, \cite[Chapter 12]{cover1999elements}, explains that the \emph{relative entropy} is the rate at which deviations from the prior path law become unlikely as the number of independent particles grows. As such, \emph{the relative-entropy minimizer is the most likely} law among all other candidates. 
The relative entropy, also known as the \emph{Kullback-Leibler} divergence, is defined between $Q$ and $P$ as
\begin{align*}
    H(P,Q) :=\int {\rm d}P  \log \frac{{\rm d}P}{{\rm d} Q},
\end{align*}
when $P$ is absolutely continuous with respect to $Q$, denoted as $P \ll Q$. This holds if every set of paths that has zero probability under $Q$ also has zero probability under $P$.

Let $\rho_t$ denote the one-time marginals of a candidate law $P$. The Schr\"odinger bridge problem can then be formulated as the following optimization problem.
\begin{subequations}\label{problem:SBP-basic}
    \begin{align}
    \underset{P: P \ll Q}{{\text{Minimize}}} \quad  H(P,Q),
\end{align}
subject to
\begin{align}\label{eq:BC-SB-basic}
    \rho_0(x) = \rho_{\rm initial}(x), \quad  \rho_T(x) = \rho_{\rm final}(x).
\end{align}
\end{subequations}
This is a convex optimization problem because the relative entropy is strictly convex in $P$, 
and the constraints are linear. 

We can write the augmented Lagrangian of the problem as
\begin{align}
    \mathcal L  &= \int {\rm d}P  \log \frac{{\rm d}P}{{\rm d} Q} \nonumber\\
    &+ \int_{\mathbb R^d} \alpha_1(x) \big [\rho_0(x) - \rho_{\rm initial}(x) \big] \, {\rm d}x \nonumber \\
    &+\int_{\mathbb R^d} \alpha_2(x) \big [\rho_T(x) - \rho_{\rm final}(x) \big] \, {\rm d}x,
\end{align}
where $\alpha_1, \alpha_2$ are Lagrange multipliers.
We let
\begin{align*}
    \Lambda := \frac{{\rm d}P}{{\rm d} Q},
\end{align*}
denoting the Radon--Nikodym derivative. Then, we can express $\mathcal L$ as
\begin{align*}
      \mathcal L  &= \int  \Lambda \Big[ \log \Lambda + \alpha_1(X_0) + \alpha_2(X_T) \Big] \, {\rm d} Q \nonumber \\
    &-\int_{\mathbb R^d} \Big[\alpha_1(x) \rho_{\rm initial}(x) + \alpha_2(x) \rho_{\rm final}(x) \Big] \, {\rm d}x.
\end{align*}

The optimal $P$, denoted below as $P^\star$, must be such that the first variation of $\mathcal L$ with respect to $\Lambda$ vanishes. Because the problem is convex, this is also a sufficient condition for optimality. The first variation is computed as
\begin{align*}
    \delta \mathcal L = \int   \Big[ \log \Lambda + \alpha_1(X_0) + \alpha_2(X_T) +1 \Big]  \delta \Lambda\, {\rm d} Q.
\end{align*}
This implies the optimal $\Lambda$ is
\begin{align*}
   \Lambda^\star = \beta_1(X_0) \beta_2(X_T)
\end{align*}
with $\beta_1(x) = e^{-1 - \alpha_1(x)}$, $\beta_2(x) = e^{-\alpha_2(x)}$. Equivalently,
\begin{align}\label{eq:opt-P}
    P^\star =  \beta_1(X_0) \beta_2(X_T) Q,
\end{align}
implying that the bridge $P^\star$ reweighs the probabilities of trajectories under $Q$ based only on the endpoint positions. Additionally, as $Q$ is a Markov law, since it corresponds to an It\^o diffusion, the law $P^\star$ retains this Markov property due to the multiplicative structure in \eqref{eq:opt-P}.

From \eqref{eq:opt-P}, we deduce that the joint probability density of the initial and final positions under $P^\star$ is
\begin{align*}
   \rho_{0,T}(x,y) =  \beta_1(x) \beta_2(y) q_0(x) k(0,x,T,y).
\end{align*}
As a consequence of the Markov property, we can write the density generated by $P^\star$ as
\begin{align}
\!\!    \rho_t(z) = \int_{\mathbb R^{2d}} \beta_1(x) \beta_2(y) q_0(x) k(0,x,t,z) k(t,z,T,y) \, {\rm d} x \, {\rm d}y \nonumber \\
=  \underbrace{\int_{\mathbb R^d} \beta_1(x)  q_0(x) k(0,x,t,z)\, {\rm d}x}_{\displaystyle  =: \hat \phi_t(z)}  \underbrace{\int_{\mathbb R^d}  \beta_2(y)   k(t,z,T,y)  \, {\rm d}y}_{\displaystyle  =:  \phi_t(z)}. \label{eq:rhot-varphit-hatvarphit}
\end{align}

Observe that $\hat \phi_t$ is a result of forward-in-time propagation of an initial condition $\hat \phi_0(x) =  \beta_1(x)  q_0(x)$, while $\phi_t$ is a result of a backward-in-time propagation of the final condition $\phi_T(y) = \beta_2(y)$. Therefore, we can express the dynamics of the pair $(\hat \phi_t,\phi_t)$ as
\begin{subequations}\label{eq:schr-sys-basic}
    \begin{align}
         \partial_t \hat \phi_t + \nabla_x \cdot (b_t \hat \phi_t) &= \frac{\sigma^2}{2} \Delta_x \hat \phi_t, \label{eq:forward} \\
        \partial_t  \phi_t + \langle b_t, \nabla_x \phi_t \rangle &= -\frac{\sigma^2}{2} \Delta_x \phi_t, \label{eq:backward}
    \end{align}
    where $\langle \cdot, \cdot \rangle$ denotes the Euclidean inner product. 
    By \eqref{eq:rhot-varphit-hatvarphit} and \eqref{eq:BC-SB-basic}, we get the boundary conditions for the previous PDEs as
\begin{align}
    \hat \phi_0(x) \phi_0(x) &= \rho_{\rm initial}(x), \label{eq:schr-sys-basic-bc1} \\
    \hat \phi_T(x) \phi_T(x) &= \rho_{\rm final}(x). \label{eq:schr-sys-basic-bc2}
\end{align}
\end{subequations}

This system of linear PDEs, coupled through their boundary conditions as shown in \eqref{eq:schr-sys-basic}, is known as the \emph{Schr\"odinger system}. Because it incorporates both forward and backward contributions from $\hat \phi_t$ and $\phi_t$, the bridge $P^\star$ is time-symmetric, favoring no particular direction of time \cite{chetrite2021schrodinger}. Further, the structure $\rho_t = \phi_t \hat\phi_t$, along with the dynamics of $\hat\phi_t$ and $\phi_t$, is formally analogous to the probability density representation in quantum wave mechanics, $\rho_t = \Psi_t \overline{\Psi}_t$.

The solution to the Schr\"odinger system can be computed numerically using what is called a Fortet--Sinkhorn iteration \cite{fortet1940resolution}. 
The iteration alternates between (i) integrating \eqref{eq:backward} backward in time from a guessed terminal condition, say $\phi_T=1$, (ii) updating the initial boundary condition $\hat \phi_0$ via \eqref{eq:schr-sys-basic-bc1}, (iii) integrating \eqref{eq:forward} forward in time from $\hat \phi_0$, and (iv) updating $\phi_T$ using $\hat \phi_T$ via \eqref{eq:schr-sys-basic-bc2}. This iteration continues until convergence is achieved.
This is summarized as follows. 
\begin{align}\label{eq:sinkhorn}
\phi_0&(x)\qquad \xleftarrow{\eqref{eq:backward}} \qquad \! \phi_T(x) \nonumber \\
\tfrac{\rho_{\rm initial}(x)}{ \phi_0(x)}& \bigg\downarrow\quad\qquad\qquad\qquad \quad\bigg\uparrow\ \tfrac{\rho_{\rm final}(x)}{\hat\phi_T(x)}\\
\hat \phi&_0(x)\qquad \!\!\xrightarrow{\eqref{eq:forward}}  \qquad \! \hat \phi_T(x).\nonumber
\end{align}

\subsection{Stochastic Control Formulation}
From \eqref{eq:rhot-varphit-hatvarphit}, we have
\begin{align*}
    \partial_t \rho_t  = \phi_t \partial_t \hat \phi_t +  \hat \phi_t \partial_t \phi_t.
\end{align*}
Further by \eqref{eq:forward} and \eqref{eq:backward}, we deduce that 
\begin{align} \label{eq:FP-posterior}
       \partial_t \rho_t + \nabla_x \cdot \big((b_t + \sigma^2 \nabla_x \log \phi_t ) \rho_t \big) &= \frac{\sigma^2}{2} \Delta_x \rho_t.
\end{align}
Note that the distinction between the dynamics of $\rho_t$ and those of $q_t$ in \eqref{eq:FP} is the additional drift term $\sigma^2 \nabla_x \log \phi_t$.

If a law $P$ is a candidate alternative to the prior $Q$, then the process governed by $P$ takes the form
\begin{align}\label{eq:posterior-Ito}
   \!\!\!\!\!   &{\rm d} X_t = \big(b_t(X_t) +  u_t(X_t) \big)\,{\rm d}t + \sigma \, {\rm d}B_t, \quad X_0 \sim \rho_0(x),
\end{align}
Applying Girsanov's theorem \cite[Theorem~8.64]{oksendal2003stochastic}
to \eqref{eq:Ito} and \eqref{eq:posterior-Ito}, we have
\begin{align}
    &\!\!\! \! \frac{{\rm d}P}{{\rm d} Q}  =  \frac{\rho_0}{q_0}  \, \exp \Big(\frac{1}{\sigma} \int_0^T \!\! \langle u_t , {\rm d} B_t^P \rangle + \frac{1}{2 \sigma^2}\int_0^T \Vert u_t\Vert^2 \, {\rm d}t \Big),\label{eq:Girsanov-1} 
\end{align}
where $B_t^{P}$ is a Brownian motion under $P$.
Taking the expectation of the logarithm of \eqref{eq:Girsanov-1} with respect to $P$, we obtain
\begin{align}\label{eq:KL-KE-equiv-2}
    H(P, Q) =   \frac{1}{2\sigma^2} \int_0^T \!\!\! \int_{\mathbb R^d} \Vert u_t(x)\Vert^2 \rho_t (x) \, {\rm d} x \, {\rm d}t,
\end{align}
where the expectation of the martingale $\int_0^T \langle u_t , {\rm d} B_t^P \rangle$ is zero. Additionally \footnote{With some abuse of notation, we intend by $H(\rho_0,q_0)$ the quantity $\int_{\mathbb R^d} \rho_0(x) \log \frac{\rho_0(x)}{q_0(x)} {\rm d} x$, involving probability densities instead of the corresponding probability measures. The same convention will be used later in the text.}, $H(\rho_0,q_0) = 0$ here because $\rho_0 = q_0 = \rho_{\rm initial}$ .

The bridge $P^\star$ minimizes $H(P, Q)$; hence, the corresponding added drift must minimize the expected control energy given on the right-hand side of \eqref{eq:KL-KE-equiv-2}. 
We can now recast the Schr\"odinger bridge problem in \eqref{problem:SBP-basic} as the following stochastic optimal control problem.
\begin{subequations}
    \begin{align}
         \underset{u_\cdot, \rho_\cdot}{{\text{Minimize}}}  \quad \frac{1}{2\sigma^2} \int_0^T \int_{\mathbb R^d} \Vert u_t(x)\Vert^2 \rho_t (x)  \, {\rm d} x\, {\rm d}t,
    \end{align}
    subject to
    \begin{align}
          \partial_t \rho_t + \nabla_x \cdot \big((b_t + u_t ) \rho_t \big) &= \frac{\sigma^2}{2} \Delta_x \rho_t,\\
             \rho_0(x) = \rho_{\rm initial}(x), \quad  \rho_T(x) &= \rho_{\rm final}(x).
    \end{align}
\end{subequations}
Throughout, the dot subscript indicates the entire time-dependent trajectory over $[0, T]$; for example, $u_\cdot = (u_t)_{t \in [0, T]}$, with analogous notation for other variables.

One can verify \cite{chen2021,caluya2021wasserstein} that optimal control has a state-feedback form, given by 
\begin{align}\label{eq:feedback-non-interacting}
    u_t^\star(x) = \sigma^2 \nabla_x \log \phi_t(x), 
\end{align}
indeed recovering \eqref{eq:FP-posterior}. The function $\phi_t$ is obtained from solving the Schr\"odinger system in \eqref{eq:schr-sys-basic}.

\subsection{Formulation for Kinematic Interacting Particle Systems}
Previous studies \cite{backhoff2020mean,chen2024density,hernandez2025propagation,orland2025schrodinger,eldesoukey2026generalized,inoue2026nonlocal} have addressed the Schr\"odinger bridge formulation for \emph{kinematic} (inertialess) interacting particle models. For a particle $i$, the kinematics read
\begin{align}\label{eq:finite-interacting}
    {\rm d}X_t^{i} = \frac{1}{N} \sum_{j =1,j\neq i}^N \Gamma(X_t^i -X_t^j) \, {\rm d} t + \sigma^i \, {\rm d} B_t^i,
\end{align}
where $1 \leq i \leq N$, and $\Gamma$ is a pairwise interaction kernel that describes how an agent $j$ influences the $i$-th agent. We assume $B_t^1, \ldots, B_t^N$ are independent and for simplicity, we take $\sigma^i := \sigma>0$.

The empirical distribution of the particles is given by
\begin{align*}
   \frac{1}{N} \sum_{i=1}^N  \delta_{X_t^i}.
\end{align*}
As $N \to  \infty$, a regime known as the \emph{mean-field} limit, particles interact through average effects. In this case, the empirical measure converges to a deterministic distribution, and the limiting kinematics are given by the McKean--Vlasov system
\begin{subequations} \label{eq:infinite-interacting}
\begin{align} 
    {\rm d}X_t = \mathcal I_t[\rho_t](X_t) \, {\rm d} t + \sigma \, {\rm d} B_t,
\end{align}
where $\rho_t$ is the probability density of collective and
\begin{align}
    \mathcal I_t[\rho_t](x) = \int_{\mathbb R^d} \Gamma(x- \tilde{x}) \rho_t(\tilde{x}) \, {\rm d} \tilde{x}.
\end{align}
\end{subequations}
Because the mean-field interaction force $\mathcal I_t$ depends on the particle density over the entire space, it is called \emph{nonlocal}. 

As before, the uncontrolled dynamics in  \eqref{eq:infinite-interacting} define a law $Q$ on the sample paths that may not reconcile given endpoint marginals $\rho_{\rm initial}$ and $\rho_{\rm final}$. Notably, in contrast to classical bridges, the prior $Q$ cannot be specified independently of the evolving density $\rho_t$ generated by the underlying path law. Specifically, the one-time marginals generated by $Q$ satisfy
\begin{align}
     \partial_t q_t + \nabla_x \cdot (\mathcal I_t[\rho_t] \, q_t) = \frac{\sigma^2}{2} \Delta_x q_t,
\end{align}
where $\rho_t$ appears explicitly in the prior dynamics. This is an important nuance of the interacting formulation compared to the non-interacting case. 

The mean-field Schr\"odinger bridge problem is then formulated as the following optimization problem \cite{backhoff2020mean}.
\begin{subequations}\label{problem:SBP-interacting}
    \begin{align}
    \underset{P: P \ll Q(P)}{{\text{Minimize}}} \quad  H(P,Q(P)), \label{eq:relative-entropy-interacting}
\end{align}
subject to
\begin{align}\label{eq:BC-SB-interacting}
    \rho_0(x) = \rho_{\rm initial}(x), \quad  \rho_T(x) = \rho_{\rm final}(x),
\end{align}
where the notation $Q(P)$ emphasizes the dependence of the prior on $P$. Unlike the classical Schr\"odinger bridge problem, this optimization is generally nonconvex.
\end{subequations}

This problem has been shown to admit an equivalent stochastic optimal control formulation \cite{backhoff2020mean}:
\begin{subequations}\label{problem:SBP-interacting-control}
    \begin{align}
         \underset{u_\cdot, \rho_\cdot}{{\text{Minimize}}}  \quad \frac{1}{2\sigma^2} \int_0^T \int_{\mathbb R^d} \Vert u_t(x)\Vert^2 \rho_t (x) \, {\rm d} x \, {\rm d}t,
    \end{align}
    subject to
    \begin{align}
          \partial_t \rho_t + \nabla_x \cdot \big(( \mathcal I_t[\rho_t] + u_t ) \rho_t \big) &= \frac{\sigma^2}{2} \Delta_x \rho_t, \label{eq:MV}\\
             \rho_0(x) = \rho_{\rm initial}(x), \quad  \rho_T(x) &= \rho_{\rm final}(x).
    \end{align}
\end{subequations}
Eq.~\eqref{eq:MV} is known as the McKean--Vlasov integro PDE.
Numerical methods for solving the formulations in \eqref{problem:SBP-interacting} and \eqref{problem:SBP-interacting-control} have been developed in \cite{chen2024density,eldesoukey2026generalized}. 

Compared to the kinematic setting described above, the bridge formulation for interacting \emph{kinetic} systems remains much less explored. In these systems, inertial particles are governed by
\begin{subequations}\label{eq:general-interacting kinetics}
    \begin{align}
    {\rm d} X_t &= V_t \, {\rm d} t, \\
    {\rm d}V_t &= \mathcal I_t [\mu_t](X_t,V_t) {\rm d} t + \sigma \, {\rm d} B_t,
\end{align}
\end{subequations}
where  $(X_t, V_t) \in \mathbb R^d \times \mathbb R^d$ denotes the position-velocity tuple, and $\mu_t= \mu_t(x,v)$ denotes the joint probability density over the phase space (positions and velocities) at time $t$. Notably, \eqref{eq:general-interacting kinetics} can be viewed as a first-order McKean--Vlasov system, just as in \eqref{eq:infinite-interacting}, but on the extended state space $\mathbb R^{2d}$, with a drift $ (V_t(X_t), \mathcal I_t [\mu_t])^\top$ and a degenerate noise intensity matrix, since the stochastic disturbances only affect the velocity directly.

Fischer \cite{fischer2014form} showed that, for broad classes of stochastic interacting systems, including those described by \eqref{eq:general-interacting kinetics}, the relative entropy $H(P, Q(P))$ quantifies the deviation from the prior and also provides the stochastic-control interpretation. However, the endpoint-constrained Schr\"odinger bridge problem for interacting kinetic systems was not addressed in that work. Here, we study this kinetic mean-field bridge problem for two types of interactions: the Cucker–Smale and Morse kernels. We also consider two forms of endpoint information: prescribed phase-space distributions and prescribed spatial marginals, leaving the endpoint velocity distributions unspecified. We obtain the nonlinear, nonlocally coupled optimality systems and propose a method to solve them numerically.

\section{Bridges over Cucker--Smale kinetics}\label{sec:CuckerSmale}
In the celebrated Cucker--Smale model,  each agent aligns its momentum with the weighted average of its neighbors', giving rise to the emergent behaviors of \emph{flocking} and \emph{schooling} \cite{carrillo2010asymptotic,carrillo2010particle,vicsek2012collective,motsch2011new,Carrillo2014,agueh2011analysis,haskovec2013flocking,choi2017emergent,gong2023crowd}. The \emph{Cucker--Smale kinetics} for the $i$-th agent (particle) are expressed as
\begin{subequations}\label{eq:dynamic-2d}
\begin{align}
 \!\! \!\!  {\rm d} X_t^i &=V_t^i \, {\rm d}t, \\
  \!\!\!\!  {\rm d} V_t^i &=\frac{1}{N}\sum_{j \neq i,j=1}^N  a(\Vert X_t^i -X_t^j \Vert) (V_t^j -V_t^i) \, {\rm d}t + \sigma \, {\rm d}B_t^i, \label{CuckerSmaleVelocityDynamics}
\end{align}
where $N$ denotes the number of agents in the swarm, and
\begin{align}\label{eq:W-cucker-smale}
    a(z) &:= \frac{K}{(1+ z^2)^\gamma}, \quad  K>0,\gamma\geq 0.
\end{align}
\end{subequations}
We assume $B_t^1, \ldots, B_t^N$ are independent.

The weights $a(\Vert X_t^i -X_t^j \Vert)$, known as \emph{communication rates} \cite{Carrillo2014,carrillo2010asymptotic}, describe how an agent's velocity is most influenced by those of its neighbors, with the effects of more distant interactions progressively waning for $\gamma>0$.  
In social and behavioral modeling contexts, these weights may represent how fast an individual adapts their strategies in response to those around them, and velocity alignment represents convergence in opinions or beliefs. Although initially formulated as a deterministic system \cite{cucker2007emergent,cucker2007mathematics}, Cucker et al. \cite{cucker2008flocking} later extended the eponymous model to account for environmental uncertainty (see also \cite{cms/1243443989,ton2014flocking,erban2016cucker}).


In the mean-field limit, the Cucker--Smale interaction at the phase-space state $(x,v)$ is determined by the collective density through
\begin{align}\label{eq:cucker-smale-forcing}
    f_t[\mu_t](x,v) = \int_{\mathbb R^{2d}} a(\Vert x- \tilde{x}\Vert) \mu_t(\tilde{x}, \tilde{v})    \big(\tilde{v} - v \big) \, {\rm d}\tilde{x}\: {\rm d}\tilde{v},
\end{align}
where $\mu_t=\mu_t(x,v)$ denotes the joint probability density over position and velocity. Hence, the corresponding mean-field prior dynamics read
\begin{subequations}\label{eq:dynamics-2d-meanfield}
\begin{align}
     {\rm d} X_t &=V_t \, {\rm d}t, \\
      {\rm d} V_t &= f_t[\mu_t](X_t,V_t) \, {\rm d}t + \sigma \, {\rm d}B_t. \label{eq:velocity-2d-meanfield}
\end{align}
\end{subequations}

Herein, we aim to steer the collective from a prescribed initial phase-space density $\mu_{\rm initial}(x,v)$ to a target final density $\mu_{\rm final}(x,v)$ via an optimal control input $u_t(x,v)$. The controlled mean-field dynamics are such that
\begin{subequations}\label{eq:dynamics-2d-meanfield-controlled}
\begin{align}
 {\rm d} X_t &=V_t \, {\rm d}t, \\
     {\rm d} V_t &= f_t[\mu_t](X_t,V_t) \, {\rm d}t +u_t(X_t,V_t) \, {\rm d}t+ \sigma \, {\rm d}B_t.
\end{align}
\end{subequations}
Equivalently, the joint phase-space density satisfies the controlled Vlasov–Fokker–Planck equation
\begin{align*}
\!\!\!  \partial_t \mu_t + \big \langle v , \nabla_x \mu_t \big \rangle   +\nabla_v  \cdot  \big(\, \mu_t  (f_t[\mu_t] &+  u_t) \big)  =  \frac{\sigma^2}{2}\Delta_v \mu_t,
\end{align*}
which is the kinetic counterpart of the McKean--Vlasov equation in \eqref{eq:MV}.

The dynamics in \eqref{eq:dynamics-2d-meanfield} define a prior evolution, while the optimal control $u_t(x,v)$ acts as a corrective drift.
From an inference viewpoint, the minimum-energy $u_t$ is the smallest required corrective drift so as to align the prior with the observed initial and final marginals, $\mu_{\rm initial}$ and $\mu_{\rm final}$. We also tackle the case where the available information on the collective, or equivalently the constraints, is the initial and final probability densities on $x$, $\rho_{\rm initial}(x)$ and $ \rho_{\rm final}(x)$ (as opposed to phase-space distributions).

\subsection{Phase-space Endpoint Constraints}\label{sec:cucker-smale-full}

We consider the following nonconvex optimization problem. 
\begin{subequations}\label{problem:cucker-smale}
\begin{align}
    \underset{u_\cdot, \mu_\cdot}{{\text{Minimize}}}   \ \ \frac{1}{2 \sigma^2} \int_0^T \!\!\!  \int_{\mathbb R^{2d}}   \Vert  u_t(x,v) \Vert^2 \mu_t(x,v) \, {\rm d} x\: {\rm d} v\: {\rm d}t, 
    \end{align} subject to
    \begin{align}
  \!\!    \partial_t \mu_t + \big \langle v , \nabla_x \mu_t \big \rangle   +\nabla_v  \cdot  \Big(\mu_t  \big(f_t[\mu_t] &+ u_t\big) \Big)  =  \frac{\sigma^2}{2}\Delta_v \mu_t, \label{eq:V-FP-cs-full} \\
         \mu_0(x, v) = \mu_{\rm initial}(x, v), \quad  \mu_T(x, v) &= \mu_{\rm final}(x, v),  \label{eq:BC-cs-full}
     \end{align}
\end{subequations}
where $\mu_{\text{initial}}$ and $\mu_{\text{final}}$ are given joint probability densities supported on subsets of $\mathbb R^{2d}$.

The augmented Lagrangian for this problem can be expressed as
\begin{align}\label{eq:cost-cs-full}
      \mathcal L_A &=   \int_0^T \!\!\!  \int_{\mathbb R^{2d}}  \,\bigg\{  \frac{1}{2 \sigma^2} \,  \Vert u_t \Vert^2 \mu_t + \lambda_t \Big[ \partial_t \mu_t + \big \langle v , \nabla_x \mu_t \big \rangle   \nonumber \\
     & + \nabla_v  \cdot  \big( \mu_t ( f_t[\mu_t] +  u_t)  \big)  -  \frac{\sigma^2}{2}\Delta_v \mu_t\Big] \bigg\} \, {\rm d} x\: {\rm d} v\: {\rm d}t,
\end{align}
where $\lambda_t= \lambda_t(x,v)$ is a Lagrange multiplier that enforces the dynamical constraint \eqref{eq:V-FP-cs-full}, i.e, a costate. 
Using integration by parts, we express \eqref{eq:cost-cs-full} as 
\begin{align} \label{eq:L-by parts}
     \mathcal L_A =   &\int_0^T \!\!\!  \int_{\mathbb R^{2d}} \bigg\{\frac{1}{2 \sigma^2} \Vert u_t \Vert^2    -\partial_t \lambda_t - \big \langle v , \nabla_x \lambda_t \big \rangle  \nonumber \\
     &  - \big \langle f_t[\mu_t] +u_t, \nabla_v \lambda_t \big \rangle  -  \frac{\sigma^2}{2}\Delta_v \lambda_t  \bigg\} \mu_t \, {\rm d} x\: {\rm d} v\: {\rm d}t \nonumber \\
     &+  \int_{\mathbb R^{2d}}  \Big(\lambda_T  \mu_T -  \lambda_0 \mu_0 \Big) \, {\rm d} x\: {\rm d} v, 
\end{align}
where the boundary terms vanish when we consider $\mu_t \to 0$ sufficiently fast as $ \Vert x \Vert, \Vert v\Vert \to \infty$.

To obtain the optimality conditions for Problem \eqref{problem:cucker-smale}, we derive the first variations of $\mathcal L_A$ with respect to $u_\cdot$ and $\mu_\cdot$, and equate them to zero. The first variation with respect to $u_\cdot$ is given by
\begin{align}\label{eq:dL-du}
 \!\!\! \!  \delta \mathcal L_A (u_\cdot, \delta u_\cdot) = \int_0^T \!\!\!  \int_{\mathbb R^{2d}} \Big \langle \Big(\frac{u_t}{\sigma^2} - \nabla_v \lambda_t\Big)  \mu_t, \delta u_t  \Big \rangle  \, {\rm d} x \, {\rm d}  v  \, {\rm d}t.
\end{align}
Setting $$\delta \mathcal L_A (u_\cdot, \delta u_\cdot) =0$$ yields the necessary condition
\begin{align}
  u_t^\star(x,v) =  \sigma^2 \,\nabla_v \lambda_t(x,v). \label{eq:u-cs-full}
\end{align}

By substituting \eqref{eq:u-cs-full} into \eqref{eq:L-by parts}, we obtain
\begin{align*}
     \mathcal L_A &=   \int_0^T \!\!\!  \int_{\mathbb R^{2d}} \Big\{ -\partial_t \lambda_t - \big \langle v , \nabla_x \lambda_t \big \rangle  - \big \langle f_t[\mu_t], \nabla_v \lambda_t \big \rangle  \nonumber \\
     &  -   \frac{\sigma^2}{2} \Vert \nabla_v\lambda_t  \Vert^2-  \frac{\sigma^2}{2}\Delta_v \lambda_t \Big\} \mu_t \, {\rm d} x\: {\rm d} v\: {\rm d}t \nonumber \\
     &+  \int_{\mathbb R^{2d}}  \Big(\lambda_T  \mu_T -  \lambda_0 \mu_0 \Big) \,  {\rm d} x\: {\rm d} v.  
\end{align*}
Letting
\begin{align}
     \Upsilon^{\rm CS} :=  -  \int_0^T \!\!\!\int_{\mathbb R^{2d}} \big \langle  f_t[\mu_t], \nabla_v \lambda_t \big\rangle \mu_t \, {\rm d} x\: {\rm d} v\: {\rm d}t, \label{eq:upsilon-cs}
\end{align}
the first variation with respect to $\mu_\cdot$ reads
\begin{align} \label{eq:dLdmu}
\!\!\!    \delta \mathcal L_A (\mu_\cdot, &\delta\mu_\cdot) = \int_0^T \!\!\!  \int_{\mathbb R^{2d}}   \Big\{ -\partial_t \lambda_t - \big \langle v , \nabla_x \lambda_t \big \rangle \nonumber \\
     & -  \frac{\sigma^2}{2} \Vert \nabla_v\lambda_t  \Vert^2-  \frac{\sigma^2}{2}\Delta_v \lambda_t \Big\} \delta \mu_t \, {\rm d} x\: {\rm d} v\: {\rm d}t  \nonumber \\
     &+ \delta  \Upsilon^{\rm CS}(\mu_\cdot, \delta\mu_\cdot).
\end{align}

Direct computation using $$a(\Vert x -\tilde{x} \Vert) = a(\Vert \tilde{x} -x\Vert)$$ yields 
    \begin{align}\label{eq:intermediate-step-in-optimality}
\!\!\! \delta  \Upsilon^{\rm CS}(\mu_\cdot, \delta\mu_\cdot)  = & \int_0^T \!\!\!  \int_{\mathbb R^{2d}}  \bigg [-\big\langle f_t[\mu_t]  , \nabla_v \lambda_t \big \rangle  \nonumber\\
&\qquad\quad+ r_t[\mu_t](x,v) \bigg]  \delta\mu_t \, {\rm d} x \:{\rm d} v \, {\rm d} t,
\end{align}
where
\begin{align}\label{eq:qt}
 \!\!\!   r_t&[\mu_t](x,v) : = \nonumber \\
    & \int_{\mathbb R^{2d}} \!a(\Vert x -\tilde{x}\Vert)\mu_t(\tilde{x},\tilde{v})  \Big \langle \tilde{v} -v,   \nabla_{\tilde v} \lambda_t(\tilde{x},\tilde{v}) \Big \rangle \, {\rm d}\tilde{x}\: {\rm d}\tilde{v}.
\end{align}
By substituting \eqref{eq:intermediate-step-in-optimality} into \eqref{eq:dLdmu}, the vanishing first variation of the Lagrangian $\mathcal L_A$ with respect to $\mu_\cdot$ yields the Hamilton--Jacobi--Bellman (HJB) equation
\begin{align}\label{eq:HJB-psi-cs-full} 
\partial_t \lambda_t + \big \langle v &, \nabla_x \lambda_t \big \rangle + \big\langle f_t[\mu_t^\star] ,\nabla_v \lambda_t  \big \rangle  \nonumber \\
 &+ \frac{\sigma^2}{2 }\Vert \nabla_v\lambda_t  \Vert^2+ \frac{\sigma^2}{2}\Delta_v \lambda_t =   r_t[\mu_t^\star](x,v),
\end{align}
for arbitrary $\delta \mu_\cdot$.
Observe that endpoint densities $\mu_0$ and $\mu_T$ are prescribed by \eqref{eq:BC-cs-full}, implying that their perturbations, $\delta \mu_0$ and $\delta \mu_T$, necessarily vanish.

Now, we define the change of variables 
\begin{align}\label{eq:log-transform}
    \lambda_t(x,v) = \log \varphi_t(x,v),
\end{align}
where $\varphi_t$ is a positive function for all $t \in [0,T]$. This transformation is known as the \emph{Cole--Hopf} \cite{cole1951quasi,hopf1950partial,teter2025hopf} or \emph{Fleming logarithmic transformation} \cite{fleming2005logarithmic}. 
By \eqref{eq:HJB-psi-cs-full} and the identity
\begin{align}
    \Vert \nabla_v \log \varphi_t \Vert^2+ \Delta_v \log \varphi_t = \frac{\Delta_v \varphi_t}{\varphi_t},
\end{align}
we deduce that $\varphi_t$ in \eqref{eq:log-transform} satisfies the reaction-advection-diffusion equation
\begin{align} \label{eq:varphit-cs-full}
    \partial_t \varphi_t + \big \langle v ,  \nabla_x \varphi_t \big \rangle   &+ \big\langle f_t[\mu_t^\star] , \nabla_v\varphi_t  \big \rangle   \nonumber \\
    &+\frac{\sigma^2}{2} \Delta_v \varphi_t = r_t[\mu_t^\star](x,v) \varphi_t.
\end{align}
Using \eqref{eq:u-cs-full} and \eqref{eq:log-transform}, the optimal control can be expressed as 
\begin{align}\label{eq:control-cs-full}
    u_t^\star(x,v) = \sigma^2 \, \nabla_v \log \varphi_t(x,v).
\end{align}
Notably, this form of the optimal controller mirrors that of the classical non-interacting Schr\"odinger bridge in \eqref{eq:feedback-non-interacting}, preserving the gradient-feedback structure. A distinction lies in the fact that the gradient in \eqref{eq:u-cs-full} is taken with respect to the velocity variable, as opposed to the position variable in the kinematic case.

We also introduce
\begin{align} \label{eq:mu-psi-phi-cs-full}
    \hat \varphi_t(x,v) := \frac{\mu_t^\star(x,v)}{\varphi_t(x,v)}.
\end{align}
By \eqref{eq:control-cs-full}, \eqref{eq:V-FP-cs-full}, and \eqref{eq:varphit-cs-full}, it follows that $\hat \varphi_t$ solves
\begin{align}\label{eq:hatvarphit-cs-full}
     \partial_t \hat\varphi_t + \big \langle  v , \nabla_x \hat\varphi_t  \big \rangle  &+ \nabla_v  \cdot (\hat\varphi_t f_t[\mu_t^\star]) \nonumber \\
     &-\frac{\sigma^2}{2} \Delta_v \hat\varphi_t = -  r_t[\mu_t^\star](x,v) \hat\varphi_t.
\end{align}
Further, by \eqref{eq:qt}, \eqref{eq:log-transform}, and \eqref{eq:mu-psi-phi-cs-full}, we can express the reaction rate as
\begin{align}
\!\!\!     r_t&(x,v) = \nonumber \\
     &\int_{\mathbb R^{2d}} a(\Vert x -\tilde{x}\Vert) \hat\varphi_t(\tilde{x}, \tilde{v}) \Big \langle \tilde{v}- v , \nabla_{\tilde{v}} \varphi_t(\tilde{x},\tilde{v}) \Big \rangle  \, {\rm d} \tilde{x} \: {\rm d}\tilde{v},
\end{align}
where we have dropped $\mu_t^\star$ from the notation. 

Let us summarize the obtained optimality conditions derived for Problem~\eqref{problem:cucker-smale} as follows.

\begin{widetext}
\begin{subequations}\label{eq:Schr-cs-full}
    \begin{align}
    \mu_t^\star(x,v) = \varphi_t(x,v)  \hat\varphi_t(x,v),& \label{eq:mut-cs-full} \\
 \partial_t \varphi_t + \big \langle v ,  \nabla_x \varphi_t \big \rangle   + \big\langle f_t[\mu_t^\star] , \nabla_v\varphi_t  \big \rangle  +\frac{\sigma^2}{2} \Delta_v \varphi_t &= r_t(x,v) \varphi_t, \label{eq:backward-cs-full} \\
   \partial_t \hat\varphi_t + \big\langle  v , \nabla_x \hat\varphi_t  \big \rangle  + \nabla_v  \cdot (\hat\varphi_t f_t[\mu_t^\star] )  -\frac{\sigma^2}{2} \Delta_v \hat\varphi_t &= -  r_t(x,v) \hat\varphi_t,   \label{eq:forward-cs-full}
\end{align}
where \begin{align}
  f_t[\mu_t^\star](x,v) &= \int_{\mathbb R^{2d}} a(\Vert x- \tilde{x}\Vert) \mu_t^\star(\tilde{x}, \tilde{v})    \big(\tilde{v} - v \big) \, {\rm d}\tilde{x}\: {\rm d} \tilde{v}, \\
  r_t(x,v) &=    \int_{\mathbb R^{2d}} a(\Vert x -\tilde{x}\Vert) \hat\varphi_t(\tilde{x}, \tilde{v})  \Big \langle \tilde{v}- v , \nabla_{\tilde{v}} \varphi_t(\tilde{x},\tilde{v}) \Big \rangle  \, {\rm d} \tilde{x} \: {\rm d}\tilde{v}, \label{eq:qt-cs-alg-full}
\end{align}
and subject to 
\begin{align} 
    \varphi_0(x,v) \hat \varphi_0(x,v) &= \mu_{\rm initial}(x,v), \label{eq:IC-sys} \\
    \varphi_T(x,v) \hat \varphi_T(x,v) &= \mu_{\rm final}(x,v). \label{eq:FC-sys}
\end{align}
\end{subequations}
\end{widetext}

The system of equations \eqref{eq:Schr-cs-full} is nonlinear and coupled both nonlocally and in time.  Eqs.~\eqref{eq:backward-cs-full} and \eqref{eq:forward-cs-full} depend on $\mu_t^\star$, which itself depends on the positive pair $(\varphi_t, \hat \varphi_t)$ as in \eqref{eq:mut-cs-full}. The reaction rate $r_t(x,v)$ in these integro-PDEs depends explicitly on $\hat \varphi_t$ and $\nabla_{\tilde{v}} \varphi_t$. As a result, the system is coupled throughout the time interval $[0,T]$, and the quantities $\varphi_t, \hat \varphi_t,\mu_t^\star$ cannot be solved for independently.

Notice that \eqref{eq:backward-cs-full} is a backward-in-time integro-PDE for $\varphi_t$, while \eqref{eq:forward-cs-full} is a forward-in-time integro-PDE for $\hat \varphi_t$.
As such, we need to solve \eqref{eq:backward-cs-full} subject to a terminal boundary condition $\varphi_T$, and \eqref{eq:forward-cs-full} subject to an initial boundary condition $\hat \varphi_0$. 
Eqs.~\eqref{eq:IC-sys} and \eqref{eq:FC-sys} show that these boundary conditions are coupled via the constraints at the endpoints. This structure resembles the Schr\"odinger system in \eqref{eq:schr-sys-basic}; however, the presence of nonlinearities, as well as nonlocal and temporal coupling, makes obtaining numerical solutions more challenging.
We propose solving \eqref{eq:Schr-cs-full} using a nested iterative computation that relies on delayed evaluations to handle the aforementioned nonlinearities. Our computational approach is detailed in Appendix~\ref{sec:sinkhorn-cs-full}. Further, Eqs.~\eqref{eq:backward-cs-full} and \eqref{eq:forward-cs-full} are adjoint. Such relation simplifies the required computations for jointly solving them numerically. We expand on this aspect in Appendix~\ref{sec:adjoint}.

\subsection{Spatial Endpoint Marginal Constraints}\label{sec:Cucker-Smale-partial}
Here, we pose a variation of Problem~\eqref{problem:cucker-smale}, where endpoint constraints are given on the spatial marginals, defined by
\begin{align}
    \rho_0(x) = \int_{\mathbb R^d} \mu_0(x,v) \, {\rm d} v, \\
    \rho_T(x) = \int_{\mathbb R^d} \mu_T(x,v) \, {\rm d} v,
\end{align}
for $x \in \mathbb R^d$.
This setting presents an inference problem in which only partial information about the system at the endpoints is available. 
The studies in \cite{chen2019multi,chiarini2022entropic} addressed this scenario for the classical kinetic Schr\"odinger bridges. Here, we extend this formulation to involve nonlocal interactions.

To this end, let $\nu_0(x,v)$ be the initial probability density associated with the prior dynamics in \eqref{eq:dynamics-2d-meanfield}. Since only the spatial marginal is prescribed at the initial time, the corresponding candidate phase-space density $\mu_0(x,v)$ is not fully specified. In other words, there is freedom in choosing the initial velocity distribution as long as
\begin{align}
\int_{\mathbb R^d} \mu_0(x,v)\,{\rm d}v=\rho_{\rm initial}(x).
\end{align}

For a candidate law $P$ on the paths of \eqref{eq:dynamics-2d-meanfield-controlled}, let $Q(P)$ denote the corresponding path law of \eqref{eq:dynamics-2d-meanfield}.
Applying Girsanov's theorem to \eqref{eq:dynamics-2d-meanfield} and \eqref{eq:dynamics-2d-meanfield-controlled} yields
\begin{align*}
    H(P,Q(P)) = H(\mu_0, \nu_0)+ \frac{1}{2 \sigma^2} \int_0^T \!\!\!  \int_{\mathbb R^{2d}} \Vert  u_t \Vert^2 \mu_t \, {\rm d} x \, {\rm d} v \, {\rm d}t.
\end{align*}
The relative entropy $H(\mu_0,\nu_0)$ accounts for the deviation of $\mu_0$ from the prior $\nu_0$. Unlike the full phase-space endpoint setting, this term does not vanish and contributes to the optimization criteria.
This leads us to the following problem formulation.
\begin{subequations}\label{problem:cucker-smale-partial}
\begin{align}
  \!\!\! \!\underset{u_\cdot, \mu_\cdot}{{\text{Minimize}}}   \  H(\mu_0, \nu_0)+ \frac{1}{2 \sigma^2} \int_0^T \!\!\!  \int_{\mathbb R^{2d}} \Vert  u_t(x,v) \Vert^2 \mu_t \, {\rm d} x \, {\rm d} v \, {\rm d}t,
    \end{align}\label{eq:cost-cs-partial-}
    subject to
 \begin{align}
\!\!\!  \partial_t \mu_t + \big \langle v , \nabla_x \mu_t \big \rangle   +\nabla_v  \cdot  \Big(\mu_t  \big(f_t[\mu_t] &+  u_t\big) \Big)  =  \frac{\sigma^2}{2}\Delta_v \mu_t, \label{eq:V-FP-cs-partial} \\
  \!\!\!      \int_{\mathbb R^d} \!\! \mu_0(x, v) \, {\rm d} v = \rho_{\rm initial}(x), \quad  \int_{\mathbb R^d} \!\!\mu_T&(x, v) \,{\rm d} v = \rho_{\rm final}(x).  \label{eq:BC-cs-partial}
     \end{align}
\end{subequations}


We write the augmented Lagrangian for Problem~\eqref{problem:cucker-smale-partial} as 
\begin{align}\label{eq:cost-cs-partial}
      \mathcal L_B &=  H(\mu_0, \nu_0)+ \mathcal L_A \nonumber \\
     & + \int_{\mathbb R^d} \kappa_0(x) \bigg[  \int_{\mathbb R^d} \!\! \mu_0(x, v) {\rm d} v - \rho_{\rm initial}(x)\bigg] \, {\rm d} x \nonumber \\
     & +  \int_{\mathbb R^d} \kappa_T(x) \bigg [\int_{\mathbb R^d} \!\!\mu_T(x, v) {\rm d} v - \rho_{\rm final}(x) \bigg]  \, {\rm d} x,
\end{align}
where $\mathcal L_A$ is the same as that in \eqref{eq:cost-cs-full}, and $\kappa_0, \kappa_T$ are Lagrange multipliers that enforce the constraints in \eqref{eq:BC-cs-partial}. Since only $\rho_0, \rho_T$ are given, the corresponding phase-space densities $\mu_0,\mu_T$ are not fully fixed; thus, their perturbations may be nonzero, unlike the scenario of Problem~\eqref{problem:cucker-smale}. 


Using integration by parts, the augmented Lagrangian $\mathcal L_B$ in \eqref{eq:cost-cs-partial} can be rewritten as
\begin{align}
     \mathcal L_B &=   H(\mu_0, \nu_0) +  \eqref{eq:L-by parts} \nonumber\\
     &+  \int_{\mathbb R^d} \kappa_0\bigg[  \int_{\mathbb R^d} \!\! \mu_0 \, {\rm d} v - \rho_{\rm initial}\bigg] \, {\rm d} x \nonumber \\
     & +  \int_{\mathbb R^d} \kappa_T \bigg [\int_{\mathbb R^d} \!\!\mu_T \, {\rm d} v - \rho_{\rm final} \bigg]  \, {\rm d} x.
\end{align}
The first variations of $\mathcal L_B$ with respect to $\mu_0$ and $\mu_T$ are
\begin{align*}
 \!\!  \!\! \delta \mathcal L_B (\mu_0,\delta \mu_0) &=  \! \int_{\mathbb R^{2d}} \! \Big[\log \frac{\mu_0}{\nu_0} \!+\! 1  \\
 &+ \kappa_0(x) - \lambda_0(x,v)\Big] \delta \mu_0 \, {\rm d} x \, {\rm d} v, \\
    \delta \mathcal L_B (\mu_T,\delta \mu_T)&= \int_{\mathbb R^{2d}} \Big[\lambda_T (x,v) + \kappa_T(x) \Big] \delta \mu_T \, {\rm d} x \,{\rm d} v.
\end{align*}
The resulting optimality conditions are
\begin{subequations}
\begin{align}\label{eq:opt-BC-cs-partial}
\mu_0^\star(x,v) &=   \underbrace{e^{-1  - \kappa_0(x)}}_{\displaystyle \hat \eta_0(x)}  \nu_0(x,v) e^{\lambda_0(x,v)}  , \\
\lambda_T(x,v) &= -\kappa_T(x) \qquad \forall v \in \mathbb R^d. \label{eq:psi-T}
\end{align}
\end{subequations}
These expressions show that $\mu_0^\star$ is constructed by the positive rescaling of the prior density $\nu_0$ by $\hat \eta_0(x)$ and $\exp(\lambda_0(x,v))$, whereas $\lambda_T$ is independent of $v$. This is a considerable departure from the setting in which the full phase-space marginals are specified at the endpoints.

The derivation steps for the necessary conditions in $u_t$ and $\mu_t$ in Section~\ref{sec:cucker-smale-full} carry over without modification, resulting in the same controller structure \eqref{eq:u-cs-full} and the HJB equation \eqref{eq:HJB-psi-cs-full}. 
Thus, the transformation $\lambda_t = \log \varphi_t$ leads to identical evolution for $\varphi_t$ and the same expression for $u_t$. Defining $\hat \varphi_t(x,v) = \mu_t^\star(x,v)/\varphi_t(x,v),$
yields the same dynamics for $\hat \varphi_t$ as in \eqref{eq:hatvarphit-cs-full}. Moreover, by \eqref{eq:opt-BC-cs-partial}, we have
\begin{subequations}\label{eq:bc-cs-partial-phis}
\begin{align}
  \mu_0^\star(x,v) &=   \underbrace{\vphantom{e^{\lambda_0(x,v)}}\nu_0(x,v)\hat\eta_0(x)}_{\displaystyle \hat\varphi_0(x,v)} \underbrace{\vphantom{\nu_0(x,v)\hat\eta_0(x)} e^{\lambda_0(x,v)}}_{\displaystyle \varphi_0(x,v)} , \\
   \varphi_T(x,v) &= \underbrace{e^{-\kappa_T(x)}}_{ =: \displaystyle   \eta_T(x)}. 
\end{align}
\end{subequations}

A salient observation is that the optimal control input at the terminal time is given by
\begin{align}\label{eq:u-cs-partial-final-time}
    u_T^\star(x,v) = \sigma^2 \nabla_v \log \varphi_T(x,v)= 0
\end{align}
since, as established in \eqref{eq:psi-T}, $\lambda_T = \log \varphi_T$ is independent of $v$. 

The optimality system corresponding to Problem~\eqref{problem:cucker-smale-partial} shares Eqs.~\eqref{eq:mut-cs-full}-\eqref{eq:qt-cs-alg-full} with the system in \eqref{eq:Schr-cs-full} but differs in the endpoint constraints. For the case at hand, the endpoint constraints are 
\begin{subequations}
\begin{align} 
   \hat \varphi_0(x,v) &= \hat\eta_0(x) \nu_0 (x,v),  \label{eq:BCs-cs-partial-all-1} \\
\hat\eta_0(x) \int_{\mathbb R^d} \nu_0(x,v)\varphi_0(x,v)  \, {\rm d} v  &= \rho_{\rm initial}(x),\\
    \varphi_T(x,v) &=\eta_T(x), \quad \forall v \in \mathbb R^d, \\
    \eta_T(x) \int_{\mathbb R^d}  \hat \varphi_T(x,v) \, {\rm d} v &= \rho_{\rm final}(x).  \label{eq:BCs-cs-partial-all-2}
\end{align}
\end{subequations}
Our iterative scheme to obtain solutions of the necessary system is detailed in Appendix~\ref{sec:sinkhorn-cs-partial}.

\section{Bridges over Morse-Potential Driven Kinetics}\label{sec:morse-}

We now consider a stochastic inertial system with both attractive and repulsive pairwise interactions. For $N$ identical particles, the model is given by
\begin{subequations}\label{eq:Morse-kinetics}
\begin{align}
 \!\!\! {\rm d} X_t^{i} &= V_t^{i} \, {\rm d}t, \\
  \!\!\!  {\rm d} V_t^{i} &=  - \frac{1}{N} \sum_{j\neq i} \nabla_{X_t^{i}} W(\Vert X_t^{i} -X_t^j\Vert) \, {\rm d} t  + \sigma \,{\rm d}B_t^{i},
\end{align}
where $1 \leq i \leq N$ and
\begin{align}
    &W(z) :=  C_R \, e^{-z/\ell_R}   - C_A \, e^{-z/\ell_A}, \\
    &\text{with } \quad C_A,C_R, \ell_A, \ell_R > 0. \nonumber
\end{align}
\end{subequations}
 As before, we assume the $d$-dimensional Brownian motions $B_t^{1}, \ldots, B_t^N$ are independent.

Eqs.~\eqref {eq:Morse-kinetics} describe agents interacting via the \emph{Morse potential}, denoted above by $W(\cdot)$
The terms $$C_R \, e^{-z/\ell_R} \quad \text{and} \quad -C_A \, e^{-z/\ell_A}$$
represent, respectively, the repulsive and attractive components of the Morse potential,  with strength $C_R$ and $C_A$ and characteristic length $\ell_R$ and $\ell_A$. The repulsive interaction mitigates collisions among agents, while attraction prevents agents from separating from the collective. 
The specifications $$C_R > C_A \quad \text{and} \quad \ell_R < \ell_A$$ amount to short-range repulsion and long-range attraction \cite{mogilner2003mutual}. This choice of parameters is suitable for describing interactions in biological and bio-inspired swarms \cite{d2006self,carrillo2013new,bolley2011stochastic,leverentz2009asymptotic,maffettone2022continuification,buttenschon2024cells}.

In the mean-field limit, the interaction experienced at the phase-space state $(x,v)$ is determined by
\begin{align}\label{eq:Morse-force}
    g_t[\rho_t](x) = -\int_{\mathbb R^d} \nabla_x W(\Vert x - \tilde{x} \Vert) \rho_t(\tilde{x})\, {\rm d} \tilde{x},
\end{align}
with 
\begin{align}
    \rho_t(x) = \int_{\mathbb R^d} \mu_t (x,v) \,  {\rm d} v.
    \end{align}
Namely, the nonlocal force generated by the Morse potential $g_t[\rho_t]$ depends exclusively on the spatial marginal $\rho_t$ and the position coordinate. As such,
\begin{align*}
    \nabla_v  \cdot  \big(\mu_t \, g_t[\rho_t] \big) = \big \langle \nabla_v \mu_t, g_t[\rho_t] \big\rangle.
\end{align*}

Given initial and final phase-space densities $\mu_{\rm initial}, \mu_{\rm final}$, we formulate the inference-control problem as follows. 
\begin{subequations}\label{problem:morse}
\begin{align}
   \underset{u_\cdot, \mu_\cdot}{{\text{Minimize}}} \quad \frac{1}{2 \sigma^2} \int_0^T \!\!\!  \int_{\mathbb R^{2d}}   \Vert  u_t(x,v) \Vert^2 \mu_t \, {\rm d} x\: {\rm d} v\: {\rm d}t, 
    \end{align} subject to
    \begin{align}
      \partial_t \mu_t + \big \langle v , \nabla_x \mu_t \big \rangle   +\nabla_v  \cdot  \Big(\mu_t  \big(g_t[\rho_t] &+ u_t\big) \Big)  = \frac{\sigma^2}{2}\Delta_v \mu_t, \label{eq:V-FP-m-full} \\
         \mu_0(x, v) = \mu_{\rm initial}(x, v), \quad  \mu_T(x, v) &= \mu_{\rm final}(x, v).  \label{eq:BC-m-full}
     \end{align}
\end{subequations}

The derivation of optimality conditions follows that of Section~\ref{sec:cucker-smale-full}, except for the contribution of the forcing $g_\cdot[\rho_\cdot]$ to the first variation with respect to $\mu_\cdot$. Specifically, the optimal controller retains the form 
\begin{align}
    u_t^\star(x,v) = \sigma^2 \nabla_v \lambda_t(x,v),
\end{align} 
where $\lambda_t$ is the Lagrange multiplier that enforces the constraint \eqref{eq:V-FP-m-full}. 
The interaction-dependent quantity $\Upsilon^{\rm CS}$ in \eqref{eq:upsilon-cs} is replaced by
\begin{align}
    \Upsilon^{\rm M} :=  -  \int_0^T \!\!\!\int_{\mathbb R^{2d}} \big \langle  g_t[\rho_t] (x), \nabla_v \lambda_t \big\rangle \mu_t \, {\rm d} x\: {\rm d} v\: {\rm d}t. \label{eq:upsilon-m}
\end{align}
Direct computation using $$\nabla_{\tilde{x}} W( \Vert x -\tilde{x} \Vert) =- \nabla_x W(\Vert x - \tilde{x}\Vert)$$ yields
    \begin{align}\label{eq:intermediate-step-in-optimality2}
\!\!\! \delta  \Upsilon^{\rm M}(\mu_\cdot, \delta\mu_\cdot)  = & \int_0^T \!\!\!  \int_{\mathbb R^{2d}}  \bigg [-\big\langle g_t[\rho_t]  , \nabla_v \lambda_t \big \rangle  \nonumber\\
&\qquad\quad+ s_t[\mu_t](x) \bigg]  \delta \mu_t \, {\rm d} x \:{\rm d} v \, {\rm d} t,
\end{align}
where
\begin{align}\label{eq:rt}
 \!\!\!   s_t&[\mu_t](x) :=   \nonumber \\
 &-\int_{\mathbb R^{2d}} \! \!  \Big \langle\nabla_x W(\Vert x -\tilde{x} \Vert),   \nabla_{\tilde v} \lambda_t(\tilde{x},\tilde{v}) \Big \rangle \mu_t(\tilde{x},\tilde{v}) \, {\rm d}\tilde{x}\: {\rm d}\tilde{v}.
\end{align}

Consequently, optimizing over $\mu_\cdot$ yields the HJB equation
\begin{align}\label{eq:HJB-psi-m-full} 
\partial_t \lambda_t + \big \langle v &, \nabla_x \lambda_t \big \rangle + \big \langle g_t[\rho_t^\star] ,\nabla_v \lambda_t  \big \rangle  \nonumber \\
 &+ \frac{\sigma^2}{2 }\Vert \nabla_v\lambda_t  \Vert^2+ \frac{\sigma^2}{2}\Delta_v \lambda_t =   s_t[\mu_t^\star](x),
\end{align}
where
\begin{align}
    \rho_t^\star(x) = \int_{\mathbb R^d} \mu_t^\star(x,v) \, {\rm d} v.
\end{align}
A similar procedure to that in Appendix~\ref{sec:sinkhorn-cs-full} yields the forward-backward optimality system for Problem~\eqref{problem:morse}.
For clarity, we state the resulting system in full.
\begin{widetext}
\begin{subequations}\label{eq:Schr-m-full}
    \begin{align}
    \mu_t^\star(x,v) = \varphi_t(x,v)  \hat\varphi_t(x,v),& \label{eq:mut-m-full} \\
 \partial_t \varphi_t + \big \langle v ,  \nabla_x \varphi_t \big \rangle   + \big \langle g_t[\rho_t^\star] , \nabla_v\varphi_t  \big \rangle  +\frac{\sigma^2}{2} \Delta_v \varphi_t &= s_t(x) \varphi_t, \label{eq:backward-m-full} \\
   \partial_t \hat\varphi_t + \big \langle  v , \nabla_x \hat\varphi_t  \big \rangle  + \big \langle g_t[\rho_t^\star] , \nabla_v \hat \varphi_t  \big \rangle  -\frac{\sigma^2}{2} \Delta_v \hat\varphi_t &= -  s_t(x) \hat\varphi_t,   \label{eq:forward-m-full}
\end{align}
where \begin{align}
   g_t[\rho_t^\star](x) &= -\int_{\mathbb R^d} \nabla_x W(\Vert x - \tilde{x} \Vert) \rho_t^\star(\tilde{x})\, {\rm d} \tilde{x}, \\
  s_t(x) &:=  -\int_{\mathbb R^{2d}} \!  \Big \langle\nabla_x W(\Vert x -\tilde{x} \Vert),   \nabla_{\tilde v} \varphi_t(\tilde{x},\tilde{v}) \Big \rangle \hat \varphi_t(\tilde{x},\tilde{v}) \, {\rm d}\tilde{x}\: {\rm d}\tilde{v}, \label{eq:rt-m-full}
\end{align}
and subject to 
\begin{align} 
    \varphi_0(x,v) \hat \varphi_0(x,v) &= \mu_{\rm initial}(x,v), \label{eq:IC-m-sys} \\
    \varphi_T(x,v) \hat \varphi_T(x,v) &= \mu_{\rm final}(x,v). \label{eq:FC-m-sys}
\end{align}
\end{subequations}
\end{widetext}

\begin{figure*}[t]
 \centering
     \begin{subfigure}[b]{1\textwidth}
        \centering
        \includegraphics[trim={2.5cm 1.5cm 3.5cm 0cm},clip,width=1\linewidth]{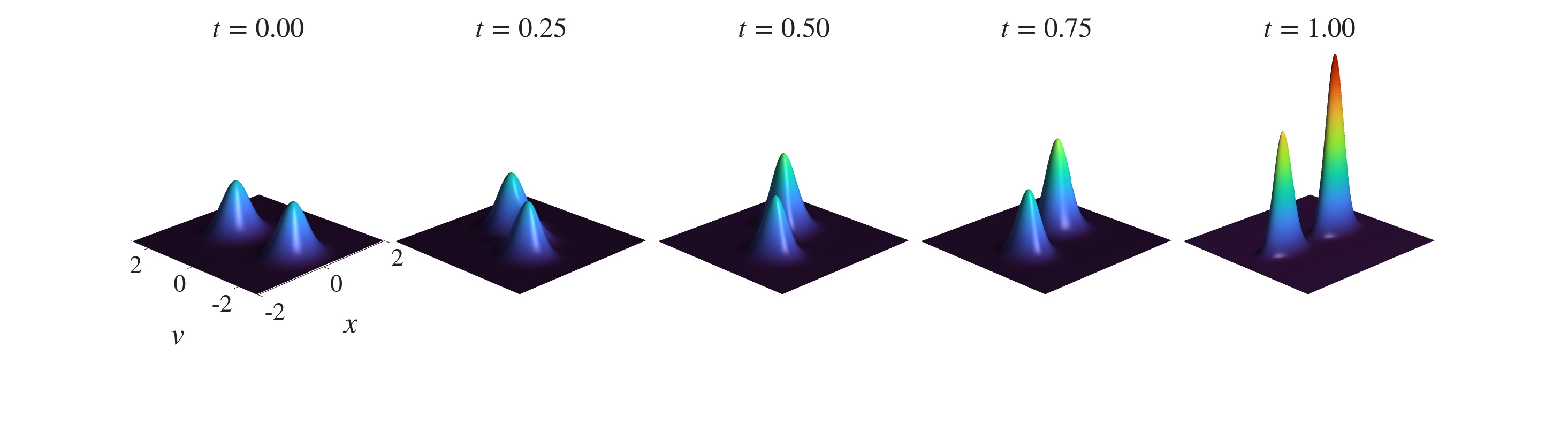}
        \includegraphics[width=1\linewidth]{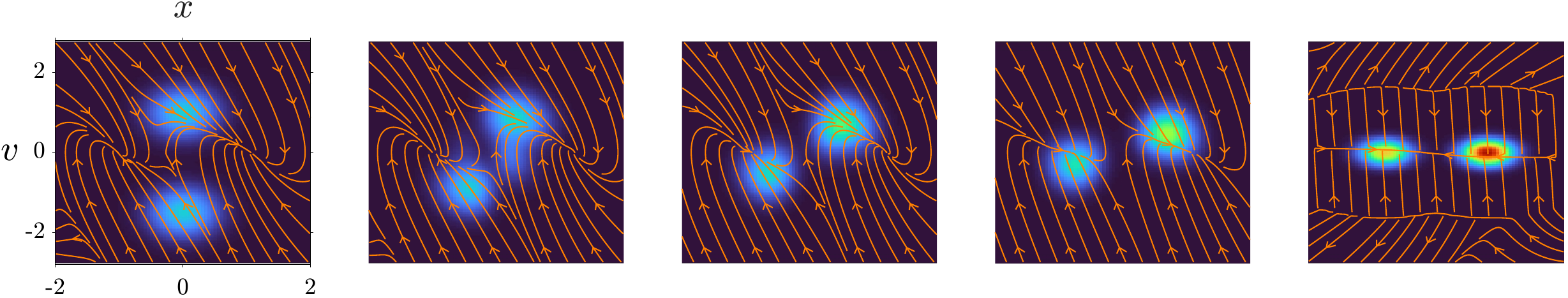}
        \includegraphics[width=\linewidth]{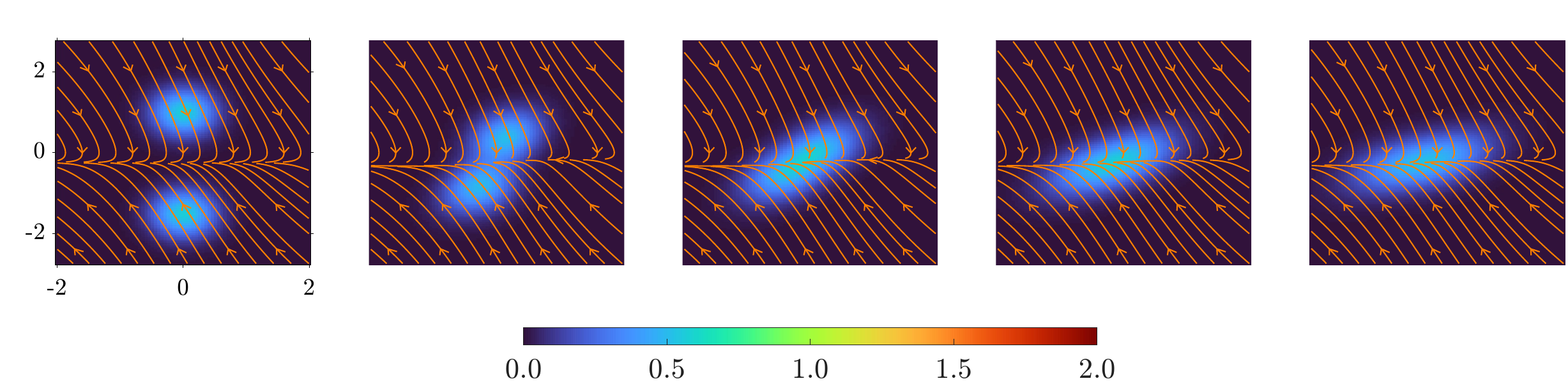}
        \centering
        \caption{At the indicated times, \emph{top row:} snapshots of the controlled phase-space density $\mu_t$, \emph{middle row:} heat maps of the controlled phase-space densities $\mu_t$, together with the instantaneous streamlines of the drift $(v, f_t[\mu_t] +u_t)$, and \emph{bottom row:} heat maps of the uncontrolled phase-space densities $\mu_t^{\rm unc}$, together with the instantaneous streamlines of the drift $(v, f_t[\mu_t^{\rm unc}])$.}
        \label{fig:Ex1-sub1}
    \end{subfigure}\\[10pt]
    \begin{subfigure}[b]{1\textwidth}
        \centering
        \includegraphics[width=\linewidth]{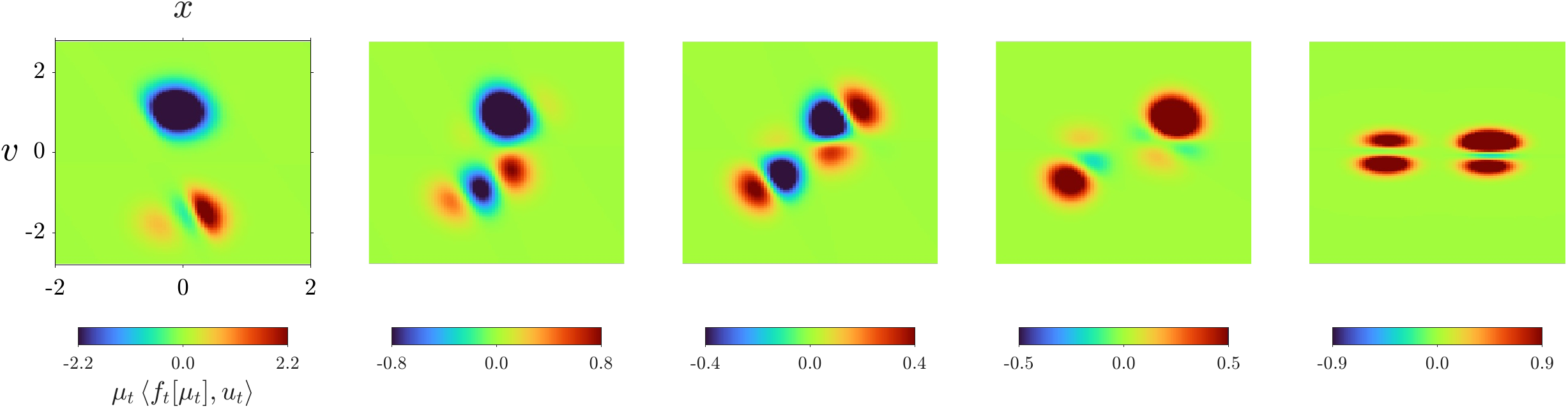}
        \caption{From left to right, heat maps of the weighted inner product between the interaction force and the controller at times $t=0,0.25,0.5,0.75$, and $1$.  The color bar limits are set according to the 97th percentile of the field's absolute values to enhance contrast, and values exceeding this range are saturated.}
        \label{fig:Ex1-sub2}
    \end{subfigure}  
    \caption{Results for the example in Section~\ref{sec:Ex-A}.}
    \label{fig:Ex1-main}
\end{figure*}

\begin{figure*}[!t]
 \centering
    \begin{subfigure}[b]{0.7\textwidth}
        \centering
         \includegraphics[trim={0cm 2cm 0cm 3.5cm},clip,width=\linewidth]{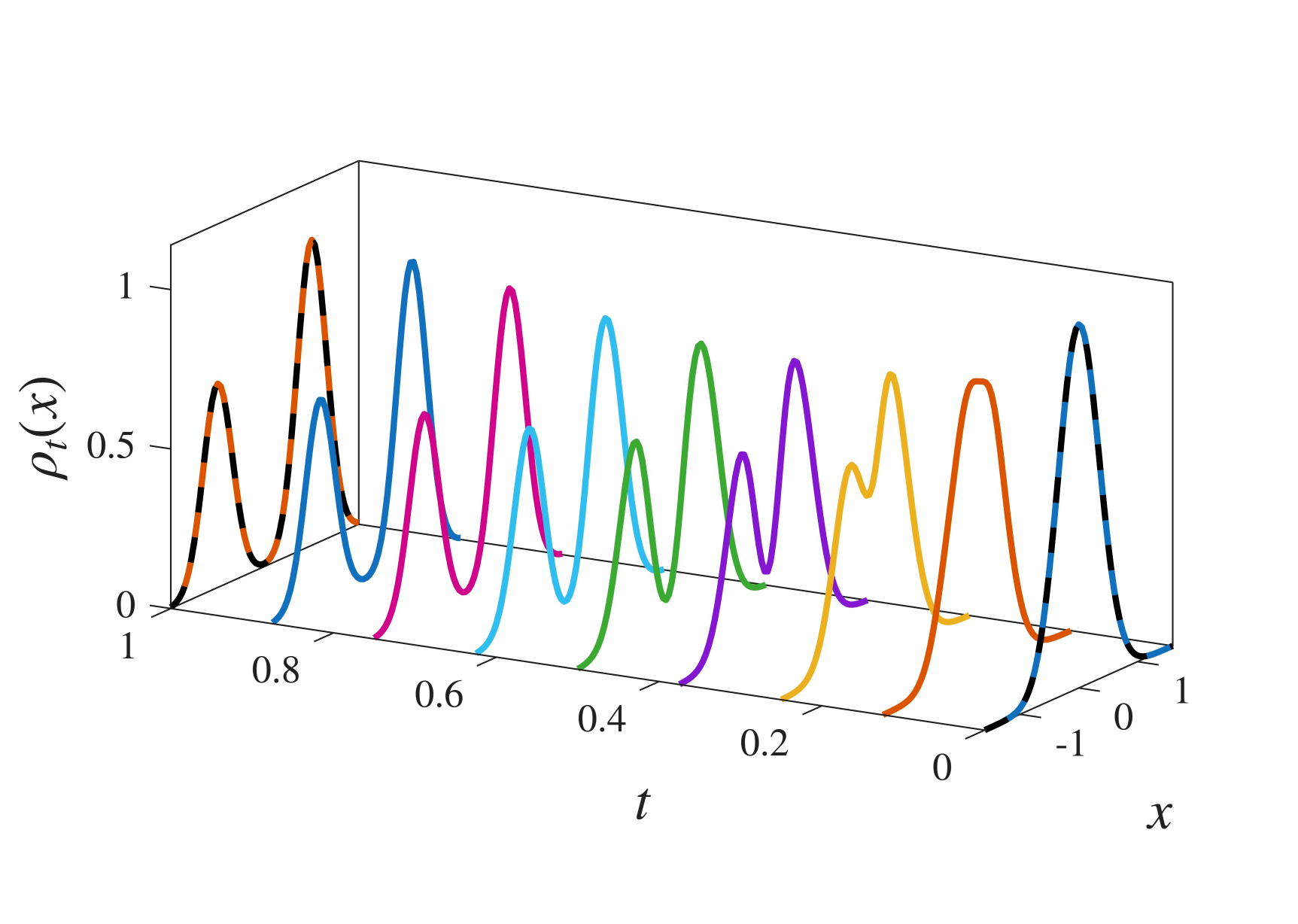}
         \vspace{2pt}
         \caption{Spatial marginals over time. Solid curves depict computed values, and dashed curves at times $t=0,1$ indicate prescribed values.} \label{fig:ex1-spatial}
    \end{subfigure}
        \begin{subfigure}[b]{0.7\textwidth}
        \centering
         \includegraphics[trim={0cm 2cm 0cm 2cm},clip,width=\linewidth]{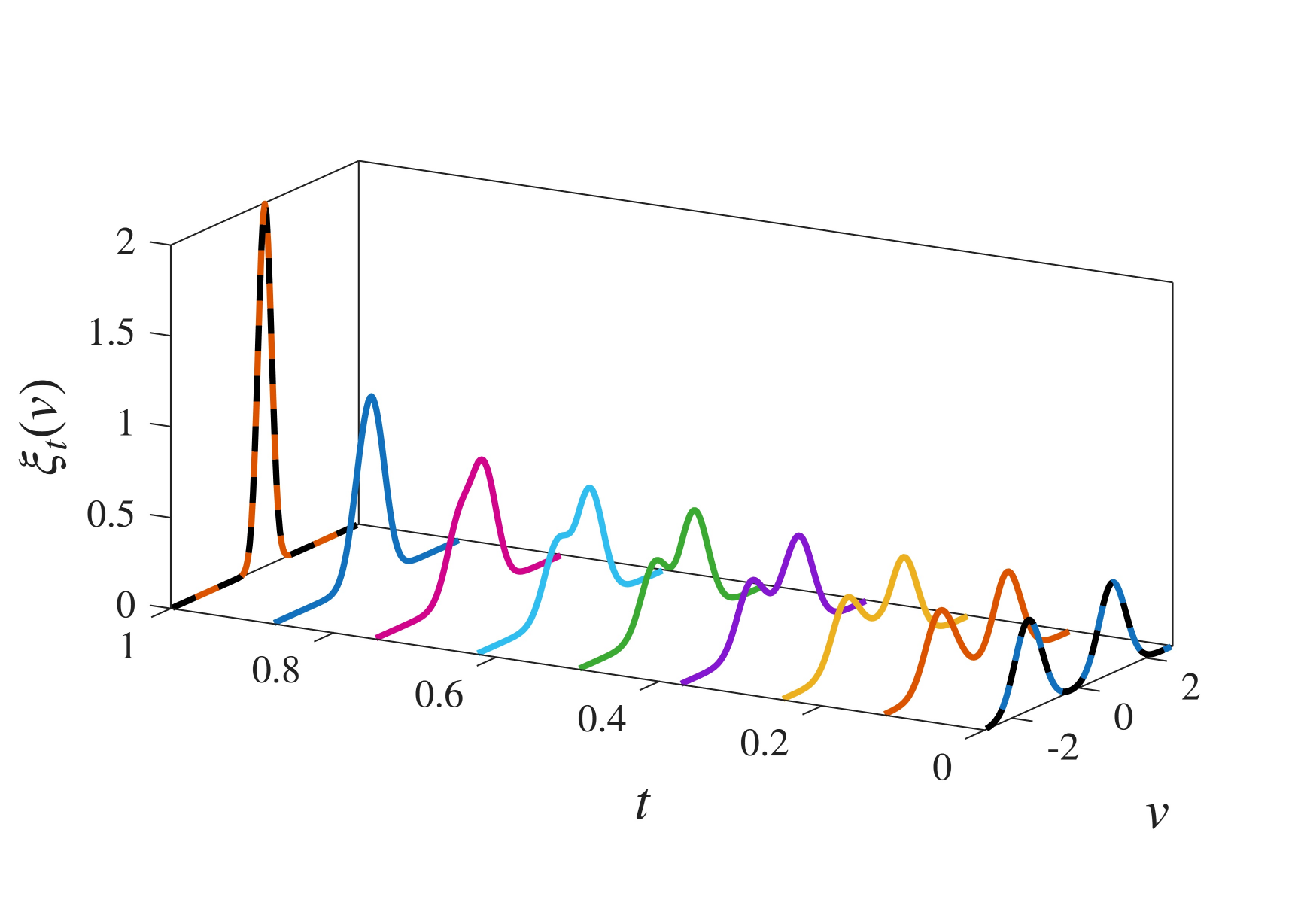}
         \vspace{2pt}
         \caption{Velocity marginals $\xi_t$ over time. Solid curves depict computed values, and dashed curves at times $t=0,1$ indicate prescribed values.}\label{fig:ex1-velocity}
    \end{subfigure}
    \vspace{15pt}\\
         \begin{subfigure}[b]{0.5\textwidth}
         \centering
         \includegraphics[width=\linewidth]{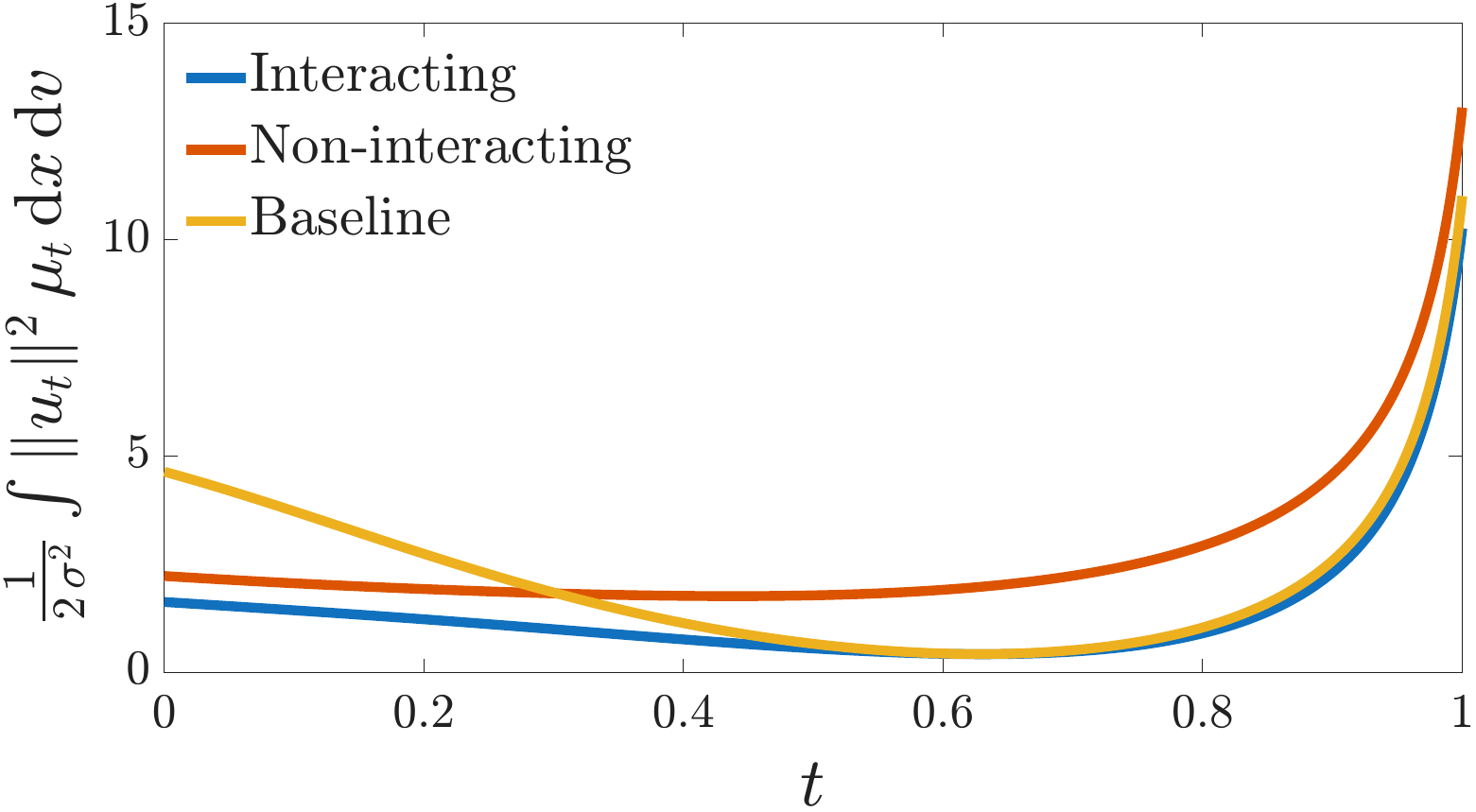}
         \vspace{2pt}
         \caption{Instantaneous cost $({1}/{2 \sigma^2})\int \Vert u_t \Vert ^2\mu_t \,{\rm d}x\:{\rm d}v$ versus time $t$ for the computed controller (interacting), the classical Schr\"odinger bridge controller (non-interacting), and the baseline controller.} \label{fig:ex1-cost}
    \end{subfigure}
    \caption{Results for the example in Section~\ref{sec:Ex-A}.}
    \label{fig:Ex1-cont}
    \end{figure*}

We next consider the inference-control problem under partial information for ensembles governed by a Morse potential. Specifically, given two marginal probability densities, $\rho_{\rm initial}$ and $\rho_{\rm final}$, both supported in $\mathbb R^d$, the corresponding optimization problem can be formulated as follows:
\begin{subequations}\label{problem:morse-partial}
\begin{align}
  \!\!\! \!\underset{u_\cdot, \mu_\cdot}{{\text{Minimize}}}   \  H(\mu_0, \nu_0)+ \frac{1}{2 \sigma^2} \int_0^T \!\!\!  \int_{\mathbb R^{2d}} \Vert  u_t(x,v) \Vert^2 \mu_t \, {\rm d} x \, {\rm d} v \, {\rm d}t, 
    \end{align}subject to \eqref{eq:V-FP-m-full} and
 \begin{align}
  \!\!\!      \int_{\mathbb R^d} \!\! \mu_0(x, v) \, {\rm d} v = \rho_{\rm initial}(x), \quad  \int_{\mathbb R^d} \!\!\mu_T&(x, v) \,{\rm d} v = \rho_{\rm final}(x).  \label{eq:BC-m-partial}
     \end{align}
\end{subequations}

Following a similar approach to that in Section~\ref{sec:Cucker-Smale-partial}, one can verify that the corresponding Schr\"odinger system for Problem~\eqref{problem:morse-partial} is described by Eqs.~\eqref{eq:mut-m-full}-\eqref{eq:rt-m-full}, together with the boundary conditions
\begin{subequations}
\begin{align} 
   \hat \varphi_0(x,v) &= \hat\eta_0(x) \nu_0 (x,v),  \label{eq:BCs-m-partial-all-1} \\
\hat\eta_0(x) \int_{\mathbb R^d} \nu_0(x,v)\varphi_0(x,v)  \, {\rm d} v  &= \rho_{\rm initial}(x),\\
    \varphi_T(x,v) &=\eta_T(x), \quad \forall v \in \mathbb R^d, \\
    \eta_T(x) \int_{\mathbb R^d}  \hat \varphi_T(x,v) \, {\rm d} v &= \rho_{\rm final}(x).  \label{eq:BCs-m-partial-all-2}
\end{align}
\end{subequations}
More details on solving the optimality systems for the Morse-type interactions are given in Appendix~\ref{sec:sinkhorn-morse}.

\section{Examples}\label{sec:Examples}
In this section, we present three numerical examples, each considering one spatial dimension and one velocity dimension. The computations are done on a periodic spatial domain and a sufficiently large truncated velocity domain so that the density near the boundaries is negligible.
\subsection{Steering the Cucker--Smale Kinetics Towards a Phase-Space Distribution}\label{sec:Ex-A}
 \begin{figure*}[t]
 \centering
     \begin{subfigure}[b]{1\textwidth}
        \centering
        \includegraphics[trim={2.5cm 1.5cm 3.5cm 0cm},clip,width=1\linewidth]{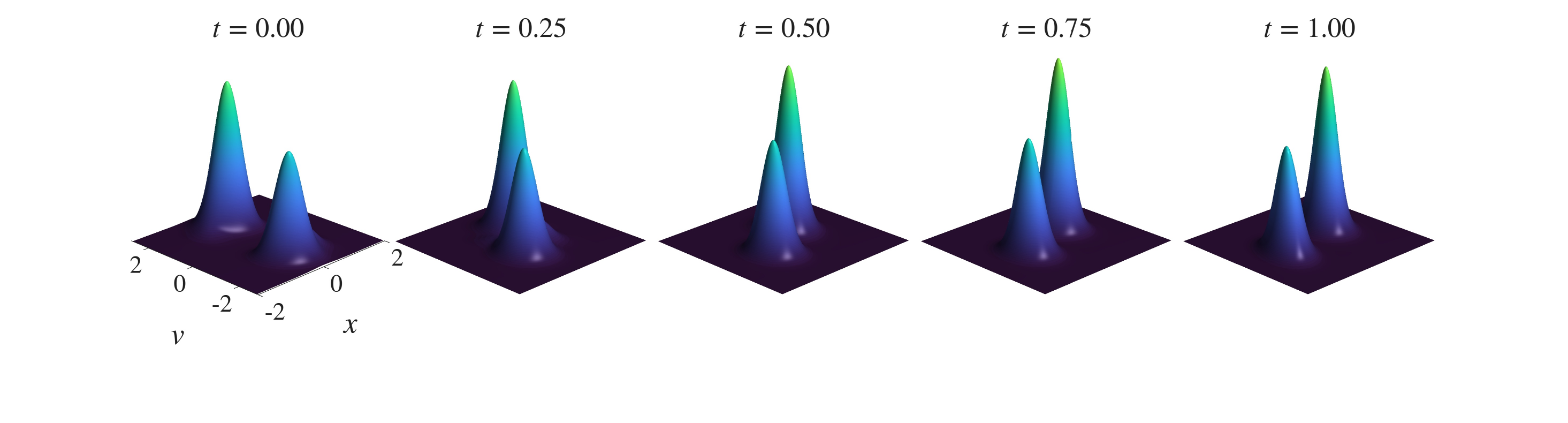}
        \includegraphics[width=1\linewidth]{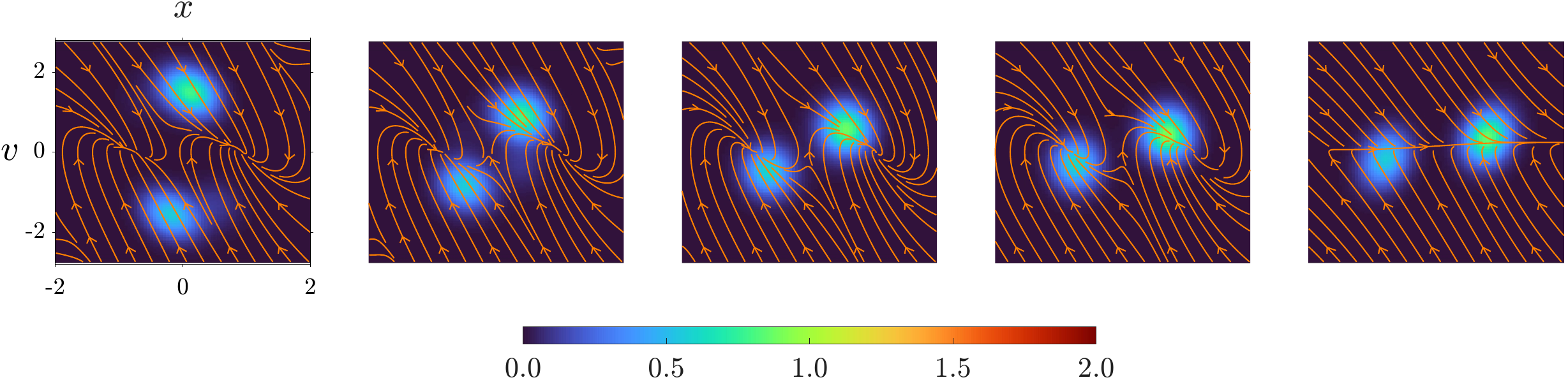}
        \centering
        \caption{At the indicated times, \emph{top row:} snapshots of the controlled phase-space density $\mu_t$ and \emph{bottom row:} heat maps of the controlled phase-space densities $\mu_t$, together with the instantaneous streamlines of the drift $(v, f_t[\mu_t] +u_t)$.}
        \label{fig:Ex2-sub1}
    \end{subfigure}\\[10pt]
    \begin{subfigure}[b]{1\textwidth}
        \centering
        \includegraphics[width=\linewidth]{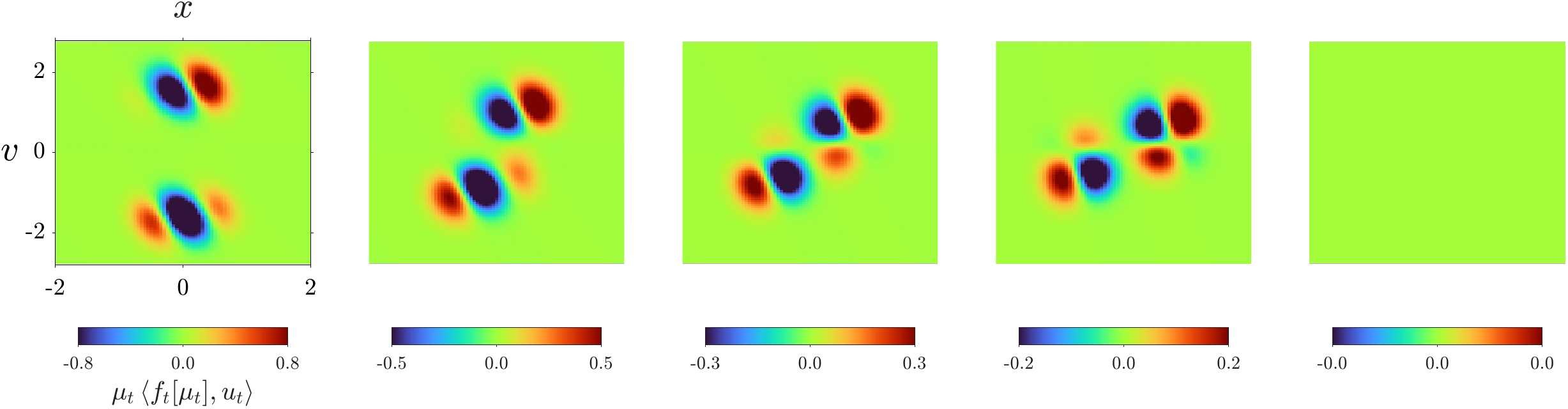}
        \caption{From left to right, heat maps of the weighted inner product between the interaction force and the controller at times $t=0,0.25,0.5,0.75$, and $1$.   The color bar limits are set according to the 97th percentile of the field's absolute values to enhance contrast, and values exceeding this range are saturated.}
        \label{fig:Ex2-sub2}
    \end{subfigure}  
    \caption{Results for the example in Section~\ref{sec:Ex-B}.}
    \label{fig:Ex2-main}
\end{figure*}

\begin{figure*}[t]
 \centering
    \begin{subfigure}[b]{0.7\textwidth}
        \centering
         \includegraphics[trim={0cm 2cm 0cm 3.5cm},clip,width=\linewidth]{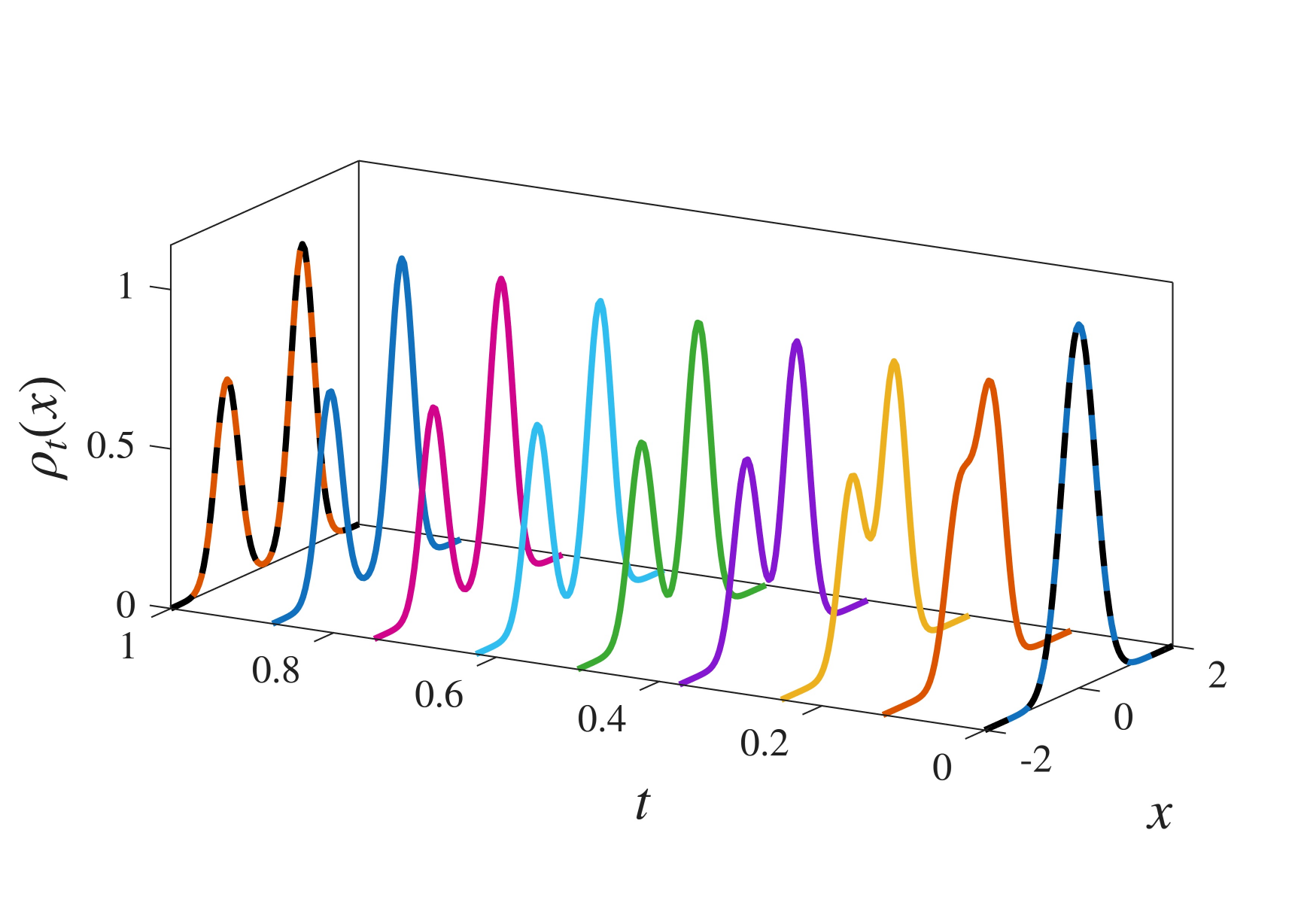}
         \vspace{2pt}
         \caption{Spatial marginals over time. Solid curves depict computed values, and dashed curves at times $t=0,1$ indicate prescribed values.}
    \end{subfigure}
        \begin{subfigure}[b]{0.7\textwidth}
        \centering
         \includegraphics[trim={0cm 2cm 0cm 2cm},clip,width=\linewidth]{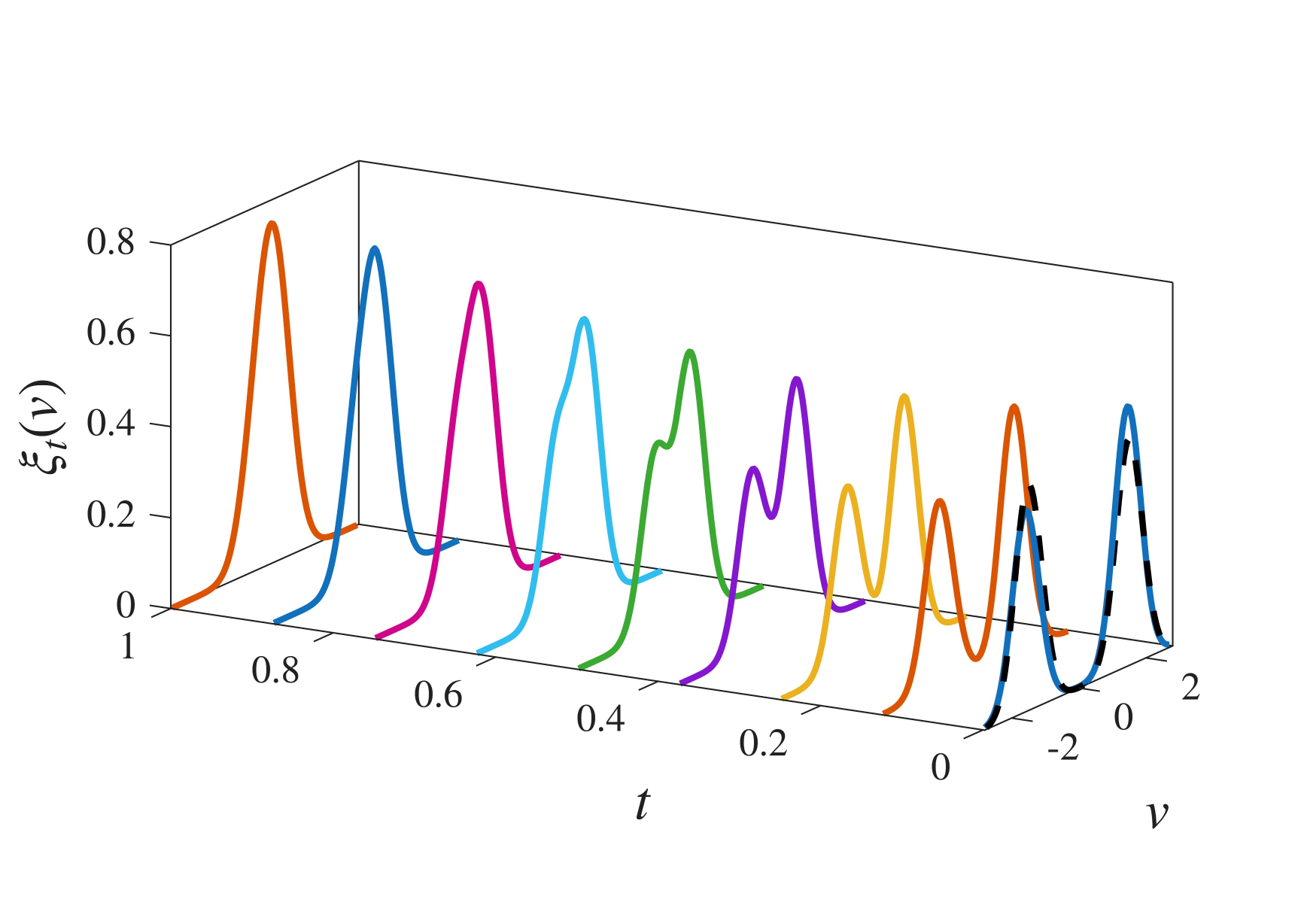}
         \vspace{2pt}
         \caption{Velocity marginals $\xi_t$ over time. Solid curves depict computed values, and the dashed curve at time $t=0$ indicates $\xi_0^{\rm prior}$.}\label{fig:ex2-velocity}
    \end{subfigure}
    \vspace{15pt}\\
         \begin{subfigure}[b]{0.5\textwidth}
         \centering
         \includegraphics[width=\linewidth]{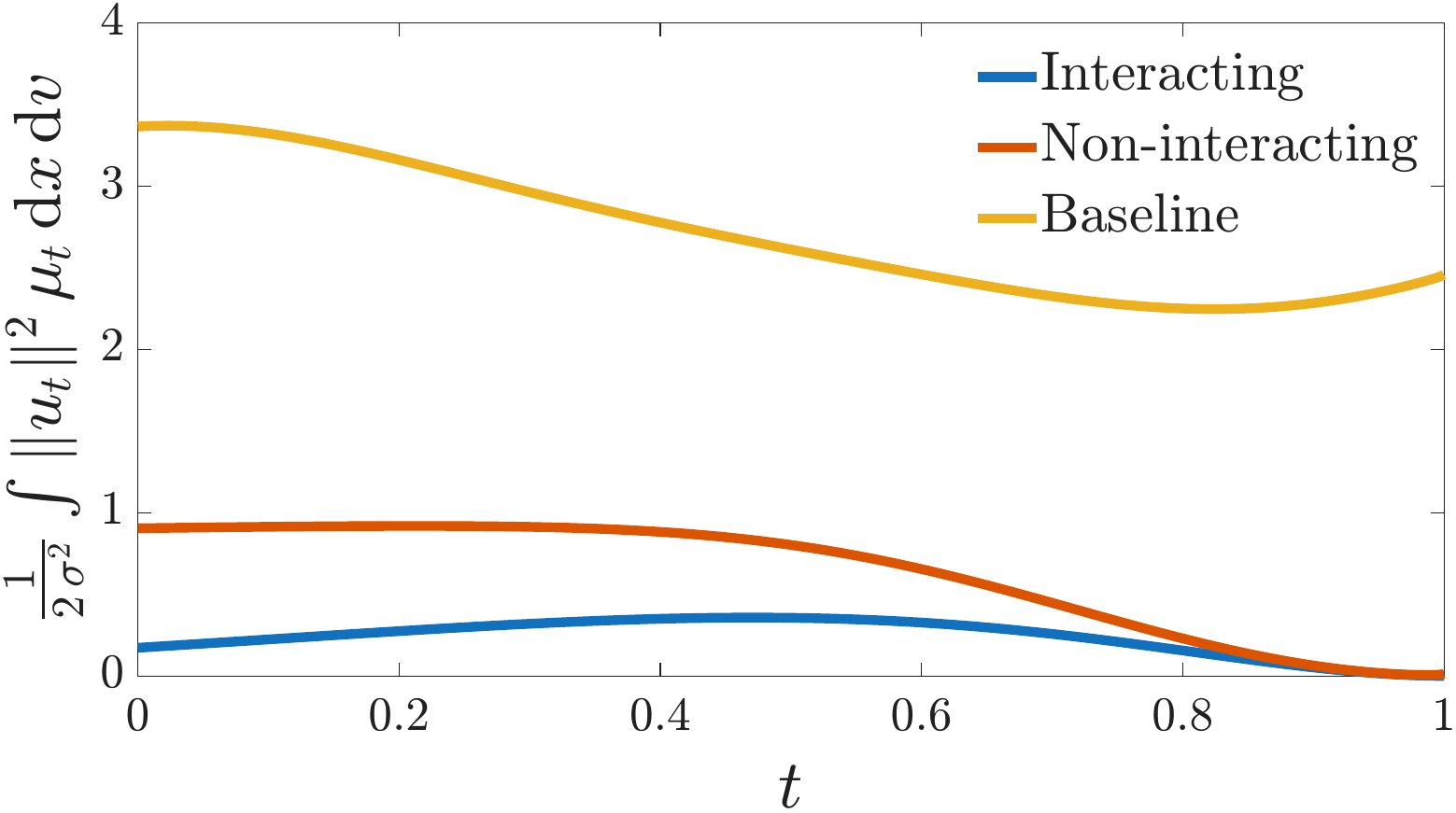}
         \vspace{2pt}
         \caption{Instantaneous cost $({1}/{2 \sigma^2})\int \Vert u_t \Vert^2\mu_t \,{\rm d}x\,{\rm d}v$ versus time $t$ for the computed controller (interacting), the classical Schr\"odinger bridge controller (non-interacting), and the baseline controller.}\label{fig:ex2-controller}
    \end{subfigure}
    \caption{Results for the example in Section~\ref{sec:Ex-B}.}
    \label{fig:Ex2-cont}
    \end{figure*}


\begin{figure*}[t]
 \centering
     \begin{subfigure}[b]{1\textwidth}
        \centering
        \includegraphics[trim={2.5cm 1.5cm 3.5cm 0cm},clip,width=1\linewidth]{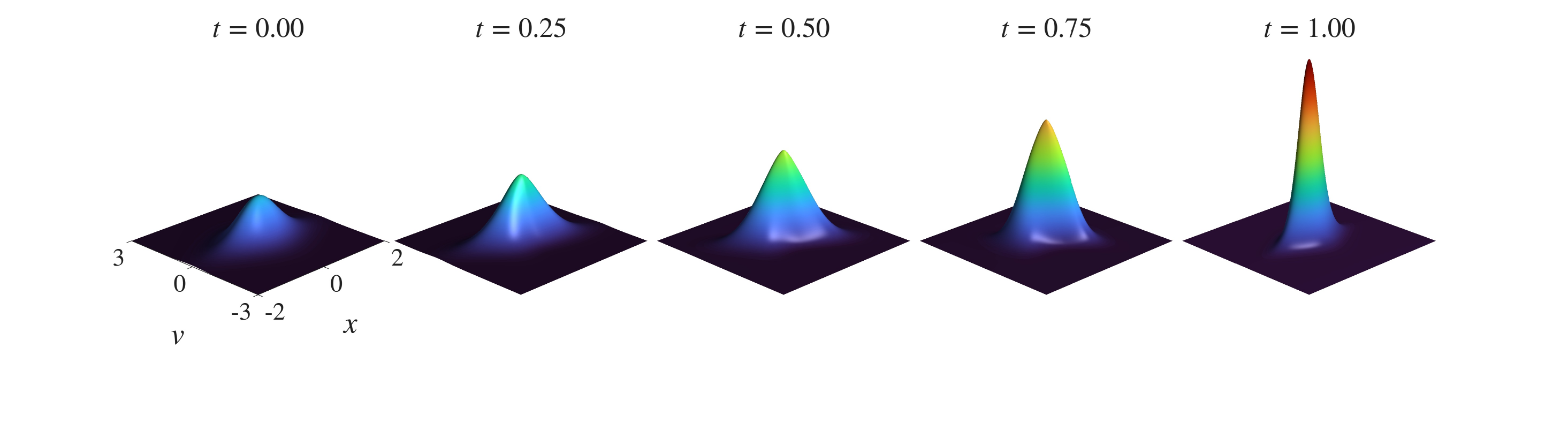}
        \centering
        \includegraphics[width=1\linewidth]{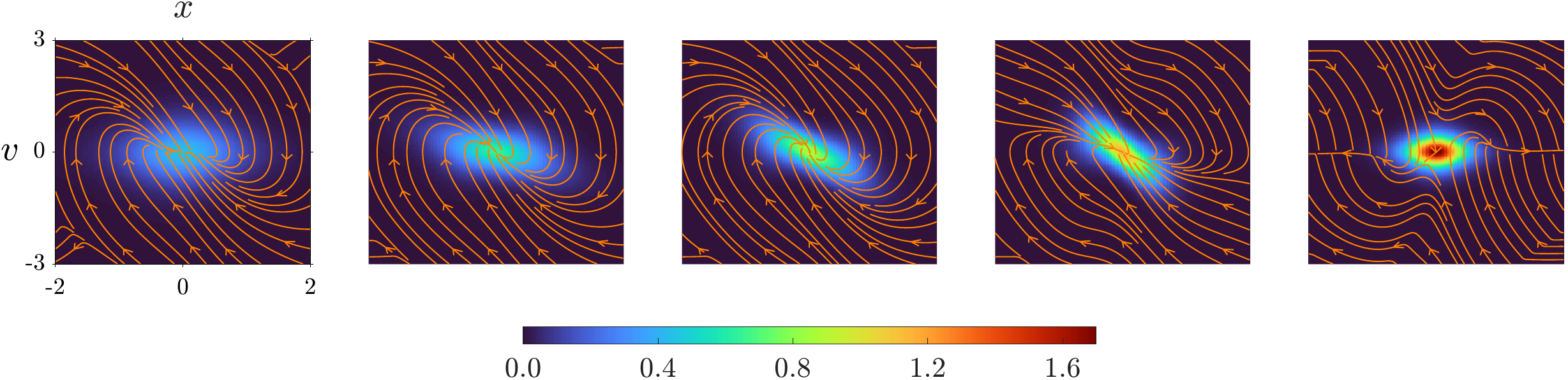}
        \caption{At the indicated times, \emph{top row:} snapshots of the controlled phase-space density $\mu_t$ and \emph{bottom row:} heat maps of the controlled phase-space densities $\mu_t$, together with the instantaneous streamlines of the drift $(v, g_t[\rho_t] +u_t)$.}
        \label{fig:Ex3-sub1}
    \end{subfigure}\\[10pt]
    \begin{subfigure}[b]{1\textwidth}
    \centering
            \includegraphics[width=\linewidth]{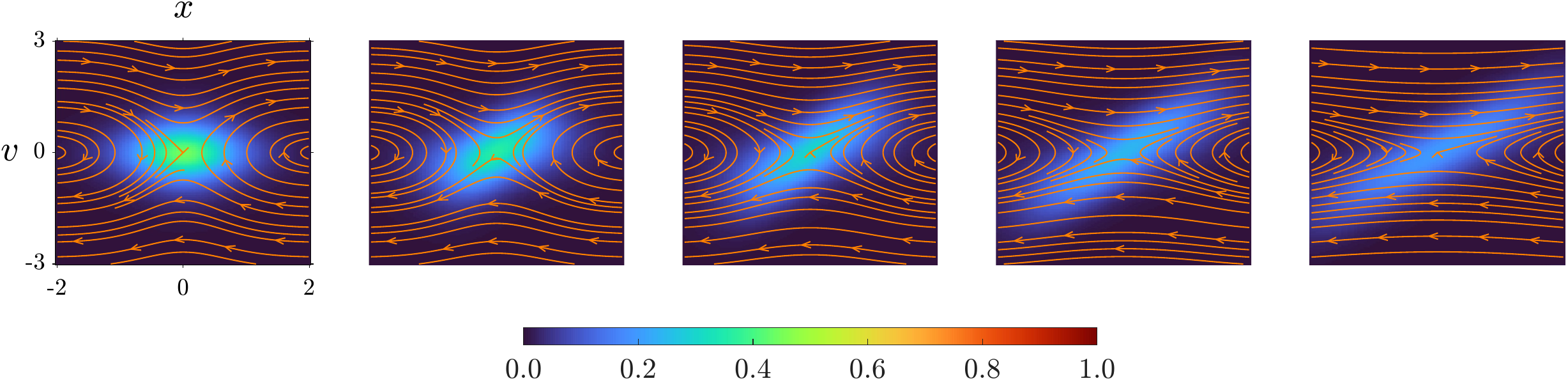}
            \caption{From left to right, heat maps of the uncontrolled phase-space densities $\mu_t^{\rm unc}$ at times $t=0,0.25,0.5,0.75$ and $1$, together with the instantaneous streamlines of the drift $(v, g_t[\rho_t^{\rm unc}])$, where $\rho_t^{\rm unc}(x)=\int \mu_t^{\rm unc}(x,v)\,{\rm d}v$ is the corresponding spatial marginal.}
            \label{fig:Ex3-sub2}
    \end{subfigure}\\[10pt]
    \begin{subfigure}[b]{1\textwidth}
        \centering
        \includegraphics[width=\linewidth]{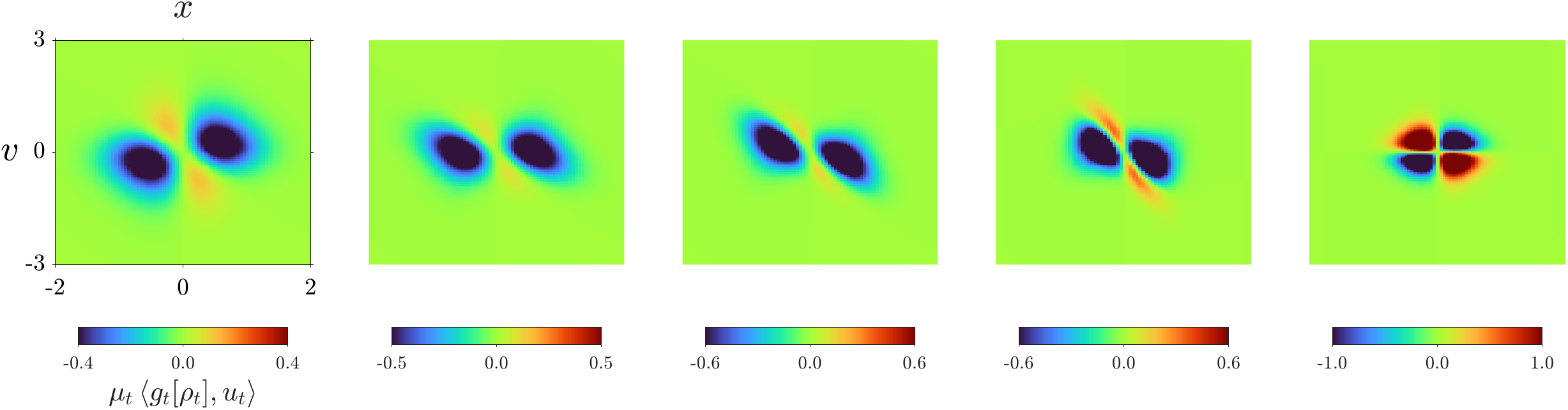}
        \caption{From left to right, heat maps of the weighted inner product of the interaction force and the controller at times $t=0,0.25,0.5,0.75$ and $1$. The color scales are set according to the 97th percentile of the field's absolute values to enhance contrast, and values exceeding this range are saturated.}
        \label{fig:Ex3-sub3}
    \end{subfigure}  
    \caption{Results for the example in Section~\ref{sec:Ex-c}.}
    \label{fig:Ex3-main}
\end{figure*}

\begin{figure*}[t]
 \centering
    \begin{subfigure}[b]{0.7\textwidth}
        \centering
         \includegraphics[trim={0cm 2cm 0cm 3.5cm},clip,width=\linewidth]{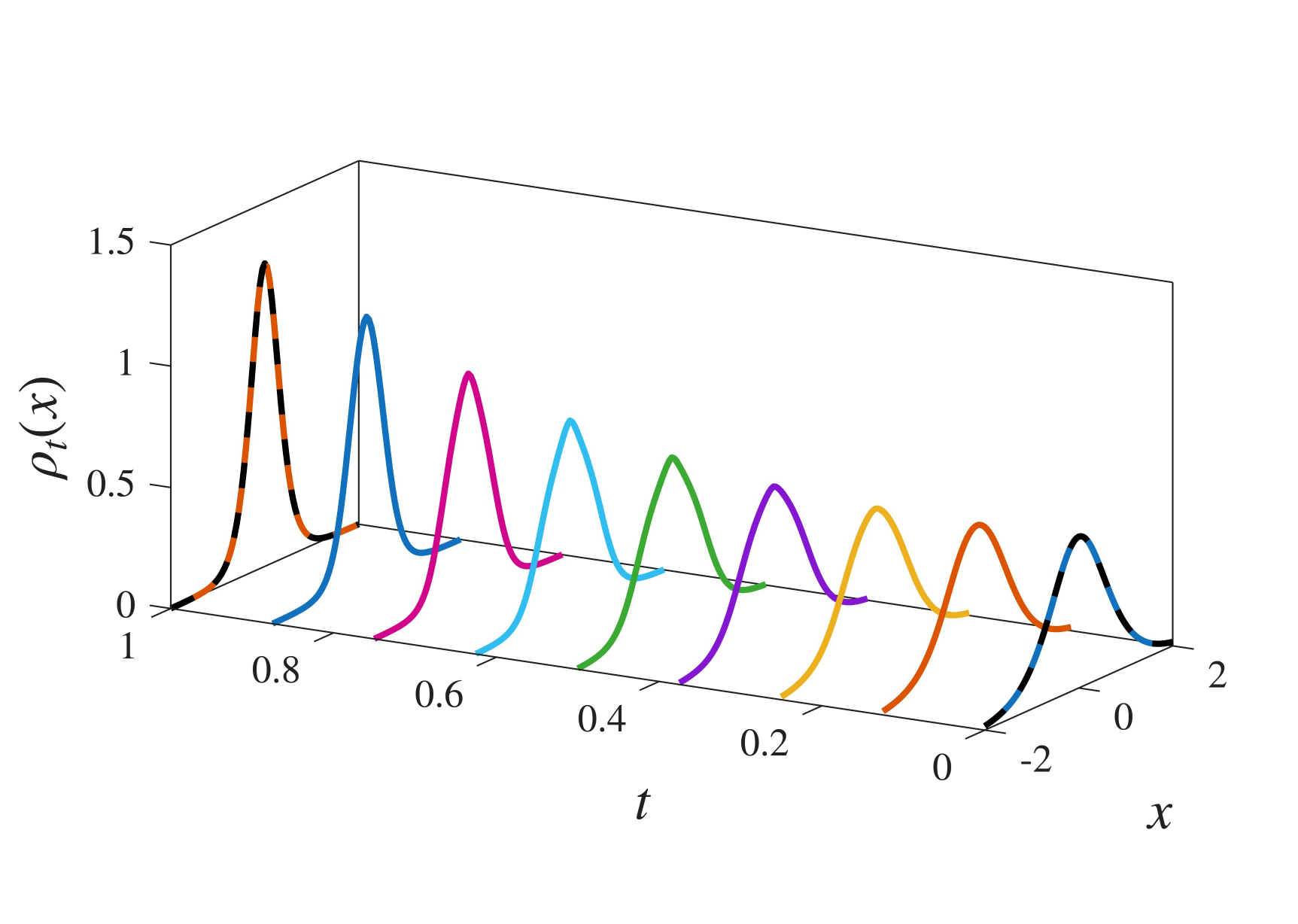}
         \vspace{2pt}
         \caption{Spatial marginals over time. Solid curves depict computed values, and dashed curves at times $t=0,1$ indicate prescribed values.} \label{fig:ex3-spatial}
    \end{subfigure}
        \begin{subfigure}[b]{0.7\textwidth}
        \centering
         \includegraphics[trim={0cm 2cm 0cm 2cm},clip,width=\linewidth]{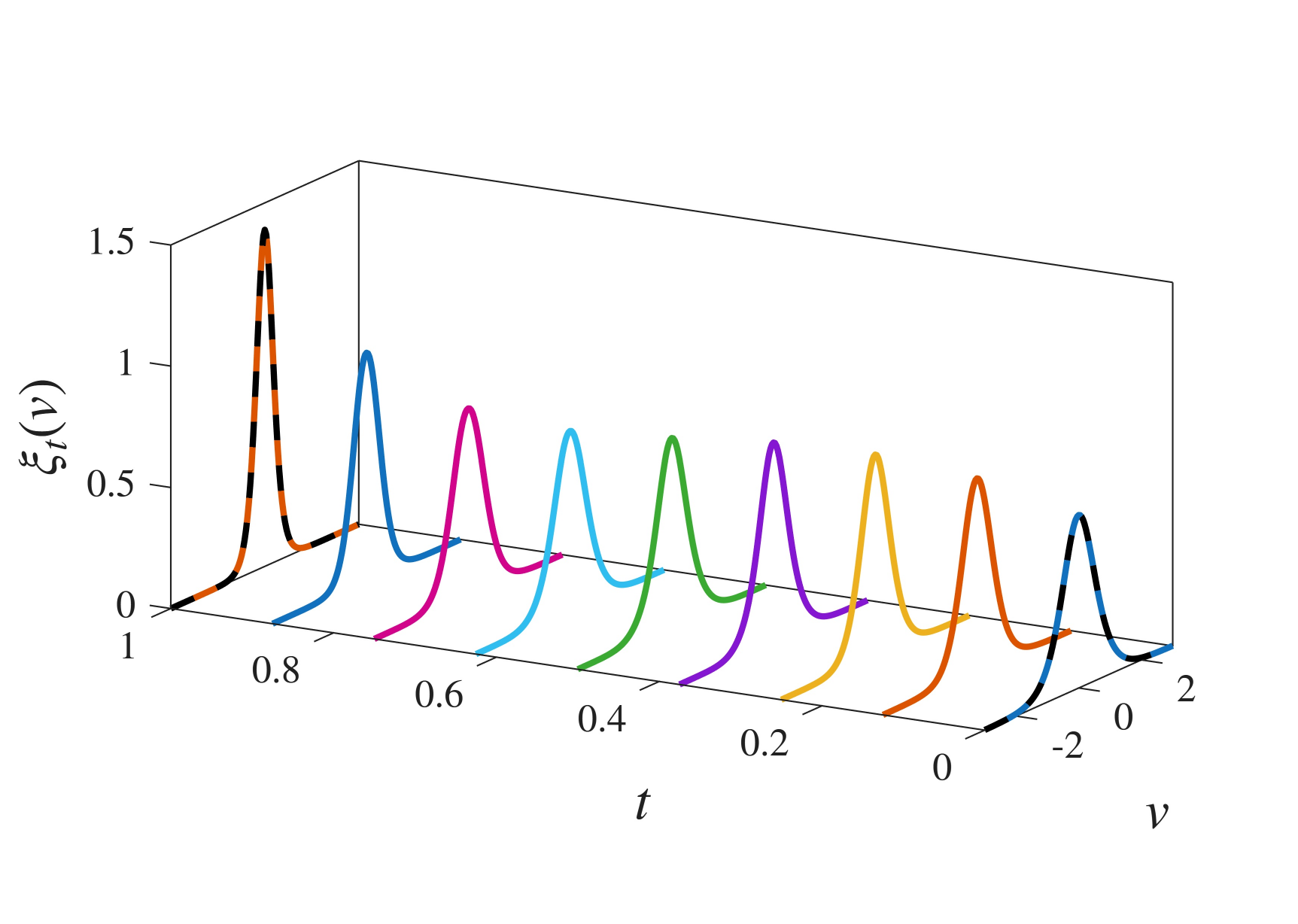} 
         \vspace{2pt}
         \caption{Velocity marginals $\xi_t$ over time. Solid curves depict computed values, and dashed curves at times $t=0,1$ indicate prescribed values.} \label{fig:ex3-velocity}
    \end{subfigure}
    \vspace{15pt}\\
         \begin{subfigure}[b]{0.5\textwidth}
         \centering
         \includegraphics[width=\linewidth]{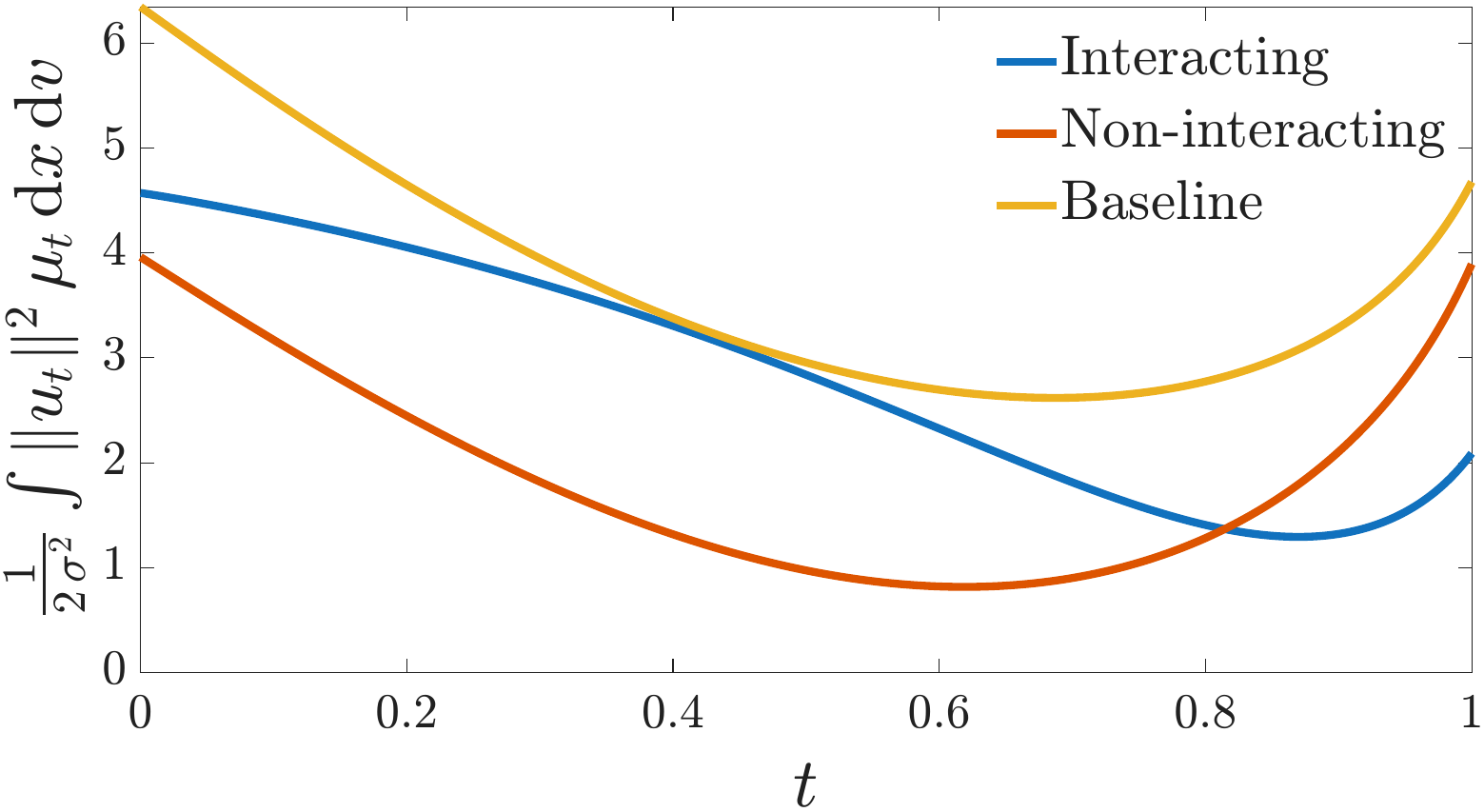}
         \vspace{2pt}
         \caption{Instantaneous cost $({1}/{2 \sigma^2})\int \Vert u_t \Vert^2\mu_t \,{\rm d}x\,{\rm d}v$ versus time $t$ for the computed controller (interacting), the classical Schr\"odinger bridge controller (non-interacting), and the baseline controller.} \label{fig:ex3-cost}
    \end{subfigure}
    \caption{Results for the example in Section~\ref{sec:Ex-c}.}
    \label{fig:Ex3-cont}
    \end{figure*}

We consider a swarm with the initial probability density given by
\begin{align}
    \mu_{\rm initial}(x,v) &= \rho_{\rm initial}(x) \xi_{\rm initial}(v), \quad \text{where } \nonumber \\
        \rho_{\rm initial}(x) &\propto e^{\frac{-x^2}{2(0.35)^2}}, \label{eq:Ex1-rho-initial} \\
     \xi_{\rm initial}(v) & \propto \Big\{e^{\frac{-(v+1.5)^2}{2(0.4)^2}} + e^{\frac{-(v-1)^2}{2(0.4)^2}}\Big\}, \nonumber
\end{align}
and $\rho_{\rm initial}$ and $\xi_{\rm initial}$ are normalized.
This distribution describes a disordered swarm with two subpopulations moving in opposite directions. More than $99.9\%$ of the velocity mass lies within the interval $[-2.8, 2.3]$, with concentrations near $-1.5$ and $1$, and a mean value of $-0.25$. We assume the swarm evolves according to the Cucker--Smale kinetics in \eqref{eq:dynamic-2d}, with $K=3$, $\gamma=0.45$, and $\sigma=1$. We consider the dynamics over $t \in [0,1]$.

The third row in FIG.~\ref{fig:Ex1-sub1} depicts the uncontrolled swarm densities $\mu_t^{\rm unc}$ that satisfy
\begin{align*}
    \partial_t \mu_t^{\rm unc}  + \langle v& , \nabla_x\mu_t^{\rm unc} \rangle  \\
    &+\nabla_v \cdot (\mu_t^{\rm unc}f_t[\mu_t^{\rm unc}] ) - \frac{ \sigma^2}{2} \Delta_v \mu_t^{\rm unc} =0,
\end{align*}
starting from $\mu_0^{\rm unc} = \mu_{\rm initial}$.
As shown in the figure, the two subpopulations gradually align their velocities over time and eventually coalesce in the phase space, demonstrating the effect of the Cucker--Smale interaction. The figure also shows streamlines that represent integral curves of instantaneous phase-space drift, $(v, f_t[\mu_t^{\rm unc}])$, where each streamline has a slope $$\frac{{\rm d}v}{{\rm d}x} = \frac{f_t[\mu_t^{\rm unc}]}{v}.$$  Generally, agents that move slower than the mean velocity are accelerated, while those moving faster than average are decelerated. The streamlines are generated using the MATLAB \texttt{streamslice} function, which prioritizes the direction and geometry of the drift vector over variations in its magnitude. Thus, streamlines at different times may look qualitatively similar, even when the interaction strength changes. 

The swarm is required to \emph{spatially split} into two groups with aligned velocities. This requirement arises, for example, when the swarm must circumvent an obstacle. The control objective is to steer the swarm toward
\begin{align}
      \mu_{\rm final}(x,v) &= \rho_{\rm final}(x) \,\xi_{\rm final}(v), \quad \text{where } \nonumber \\
    \rho_{\rm final}(x) & \propto \Big\{ 0.4 \, e^{\frac{-(x+0.8)^2}{2(0.25)^2}}  + 0.6 \, e^{\frac{-(x-0.8)^2}{2(0.25)^2}} \Big\}, \label{eq:Ex1-rho-final}\\
    \xi_{\rm final}(v) & \propto e^{- \frac{v^2}{2(0.2)^2}}, \nonumber
\end{align}
where $\rho_{\rm final},\xi_{\rm final}$ are normalized.  That is, the controller must achieve spatial splitting, with more agents on the right, while synchronizing velocities at the final time. 

With $\mu_{\rm initial}$ and $\mu_{\rm final}$ prescribed, we solve the system of equations \eqref{eq:Schr-cs-full} numerically using the method proposed in Appendix~\ref{sec:sinkhorn-cs-full} with $\theta =1$ and convergence tolerance $\epsilon=10^{-5}$. The results are presented in FIG.~\ref{fig:Ex1-main} and FIG.~\ref{fig:Ex1-cont}. The top two rows in FIG.~\ref{fig:Ex1-sub1} depict the evolution of phase-space distributions. Streamlines in the second row vary in space and time under the influence of control, showing regulated deceleration of the group moving to the right (positive $v$) and acceleration of the one moving to the left (negative $v$). 

FIG.~\ref{fig:Ex1-sub2} shows the control effect via the weighted inner product of the controller and the Cucker--Smale forcing. Initially, the controller opposes the Cucker--Smale interaction, accelerating the right subgroup and part of the left, which later merges with the right group. This acceleration facilitates spatial splitting by giving the associated masses greater rightward velocity. After the groups are nearly split, the controller, assisted by the interaction, brakes the swarm.

FIG.~\ref{fig:ex1-spatial} depicts the evolution of spatial marginals and FIG.~\ref{fig:ex1-velocity} the evolution of velocity marginals, showing that the computed marginals match the prescribed ones. FIG.~\ref{fig:ex1-cost} compares the instantaneous cost of the computed interacting controller with that with that of the controller $u_t^{\rm NSB}$, associated with the classical non-interacting second-order Schr\"odinger bridge, corresponding to $a=0$, where $$u_t^{\rm NSB}(x,v) = \sigma^2 \nabla_v \log \phi_t(x,v).$$ Here we use the symbol $u_t^{\rm NSB}$ to distinguish this controller from the other controllers. The computation of $\phi_t$ is detailed in Appendix~\ref{sec:sinkhorn-cs-full}.
The figure also shows the instantaneous cost of a basic control $$u_t^{\rm baseline} = -f_t[\mu_t^{\rm NSB}]+u_t^{\rm NSB},$$ which cancels the interaction and uses the non-interacting solution. Both the computed controller that corresponds to the system \eqref{eq:Schr-cs-full} and $u_t^{\rm baseline}$ are admissible controls that steer the interacting system between the given marginals, while $u_t^{\rm NSB}$ is the non-interacting benchmark. The accumulated cost values for the computed interacting, non-interacting, and baseline controllers are $1.34$, $2.69$, and $2.07$, respectively. Among the three, our proposed controller has the lowest accumulated cost since it utilizes the interaction force when advantageous.

\subsection{Inference with the Cucker--Smale Kinetics and Partial Information}\label{sec:Ex-B}

We consider a variant of the previous example, aiming to determine a drift vector field that steers the Cucker--Smale dynamics from $\rho_{\rm initial}$ to $\rho_{\rm final}$ given in \eqref{eq:Ex1-rho-initial} and \eqref{eq:Ex1-rho-final}. As a prior, we set
\begin{align*}
    \nu_0(x,v) &= \rho_{\rm initial}(x) \,\xi_0^{\rm prior}(v), \quad \text{where} \\
     \xi_0^{\rm prior}(v) & \propto \Big\{ e^{\frac{-(v+1.5)^2}{2(0.4)^2}} + e^{\frac{-(v-1.5)^2}{2(0.4)^2}} \Big\},
\end{align*}
where $\xi_0^{\rm prior}$ is normalized. This prior represents our belief that particles initially move right or left with equal probability. Our goal is to update this belief and obtain an evolution consistent with the spatial marginals.

With all other parameters unchanged from the previous example, we solve the optimality system for Problem~\eqref{problem:cucker-smale-partial} following the method in Appendix~\ref{sec:sinkhorn-cs-partial} with convergence tolerance $\epsilon=10^ {-5}$.
Results are presented in FIGs.~\ref{fig:Ex2-main} and \ref{fig:Ex2-cont}. 

As seen in FIG.~\ref{fig:Ex2-sub1} and the velocity marginals in FIG.~\ref{fig:ex2-velocity}, the final velocity distribution has nonzero mean and is less concentrated compared to that in the previous example (cf. FIGs.~\ref{fig:Ex1-main}, \ref{fig:Ex1-cont}). This is because there is no constraint on velocities at the final time.   

FIG.~\ref{fig:Ex2-sub2} indicates that the obtained drift (controller) accelerates the mass farthest from the target states $x=-0.8$ and $x=0.8$. Its effect is nearly symmetric for $v>0$ and $v<0$. At the final time, the control vanishes, which is consistent with \eqref{eq:u-cs-partial-final-time}. FIG.~\ref{fig:ex2-controller} includes comparison with other controllers defined analogously to those in the previous example. The accumulated cost values for the computed interacting, non-interacting, and baseline controllers are $0.24$, $0.63$, and $2.71$, respectively. Once again, the controller consistent with the necessary optimality system has the lowest cost among the three. 

\subsection{Steering the Morse-potential-driven Dynamics Towards a Target Phase-Space Distribution} \label{sec:Ex-c}
We consider a population whose initial distribution is expressed as a product of logistic distributions:
\begin{align*}
    \mu_{\rm initial}(x,v) &= \rho_{\rm initial}(x) \, \xi_{\rm initial}(v), \quad \text{where } \\
    \rho_{\rm initial}(x) &=  \frac{1}{1.6} \, {\rm sech}^2 \Big(\frac{x}{0.8}\Big), \\
    \xi_{\rm initial}(v) &=  \frac{1}{1.4} \,{\rm sech}^2\left(\frac{v}{0.7}\right).
\end{align*}
We assume that this population follows the dynamics in \eqref{eq:Morse-kinetics} with $\sigma = 1/\sqrt{2}$ and subject to the short-range repulsion and long-range attraction defined by
\begin{align*}
    C_R = 5, \quad C_A = 1.3, \quad \ell_R = 0.4, \quad \ell_A =1. 
\end{align*}
FIG.~\ref{fig:Ex3-sub2} illustrates the evolution of the uncontrolled population, showing how the repulsive component of the Morse interaction disperses the initially very concentrated group. 

We aim to steer this population towards the final distribution:
\begin{align}
     \mu_{\rm final}(x,v) &= \rho_{\rm final}(x) \xi_{\rm final}(v), \quad \text{where }\\
    \rho_{\rm final}(x) &=  \frac{1}{0.8} \, {\rm sech}^2 \Big(\frac{x}{0.4}\Big), \\
    \xi_{\rm final}(v) &=  \frac{1}{0.72} \,{\rm sech}^2 \left(\frac{v}{0.36}\right).
\end{align}
This target configuration represents a population that is more concentrated in space, with velocities more aligned, relative to the initial distribution, thereby mimicking a cohesive flock characterized by simultaneous spatial aggregation and velocity alignment. The results presented in FIGs.~\ref{fig:Ex3-sub1}, \ref{fig:Ex3-sub3}, and \ref{fig:Ex3-cont} are obtained by applying the numerical method in Appendix~\ref{sec:sinkhorn-morse} with $\theta=0.8$ and convergence tolerance $\epsilon =10^{-5}$. 

FIG.~\ref{fig:Ex3-sub3} indicates that the controller counteracts the interaction force for the most part so as to draw the population closer together. FIGs.~\ref{fig:ex3-spatial} and \ref{fig:ex3-velocity} show that the spatial and velocity marginals match the specified distributions at $t=0$ and $t=1$. Additionally, we compared the cost of the proposed controller that corresponds to \eqref{eq:Schr-m-full} vis-à-vis the costs from other controllers defined in the previous two examples. The accumulated cost values for the computed interacting, non-interacting, and baseline controllers are $2.96$, $1.94$, and $3.71$, respectively. Since the interacting control has to counteract the repulsive interaction for the most part, its cost is expected to be higher than that of the non-interacting controller, which does not need to overcome any adversarial forcing. Nevertheless, the interacting controller has a lower cost than the baseline controller.

\section{Conclusions}
In this work, we have formulated a stochastic optimal control problem for regulating interacting inertial systems in the mean-field limit to satisfy prescribed endpoint distributions. The present formulation, anchored in Schr\"odinger's paradigm of bridges, can be viewed as an inverse problem that seeks a suitable update of the interacting dynamics required to match the given endpoint distributions, regarded as observations. We considered two interaction mechanisms: Cucker--Smale alignment and Morse attraction--repulsion, and for each case, we addressed two types of endpoint specifications/information: phase-space marginals and spatial marginals alone. Interestingly, the obtained necessary optimality systems preserve the multiplicative structure of phase-space probability densities $\mu_t = \varphi_t \hat \varphi_t$, with $\varphi_t$ and $\hat \varphi_t$ satisfying backward and forward equations, respectively, albeit nonlinear and nonlocally coupled through the mean-field interactions.  

To numerically attain solutions to the derived systems, we introduced a nested fixed-point iterative approach that relies on delayed evaluations to handle the nonlinearities.
Our numerical examples, aside from attesting to the convergence of the proposed approach, demonstrate that the controller adapts its strategy based on the interaction, either opposing it when it acts against the prescribed objective or exploiting it when it becomes favorable. In the Cucker--Smale examples, where the interaction benefits the steering task, the obtained controller has costs of approximately $50\%$ and $62\%$, compared to the control strategies associated with classical Schr\"odinger bridges, which do not take interaction into account. 

The optimization problems considered are nonconvex since the interaction depends on the evolving swarm density; thus, the systems derived here provide first-order necessary conditions for optimality, and the computed solutions are stationary points rather than global minima. A future research direction is therefore to investigate second-order conditions to characterize the nature of the stationary points. Another research direction is to extend the present framework to entail time-resolved distributions (a multi-marginal formulation). This is relevant to scenarios where time series of population snapshots are available, such as in developmental trajectory inference \cite{chizat2022trajectory}. Beyond the mean-field limit, investigating finite populations as well as accounting for uncertainty in the interaction mechanisms, are also of research interest. 

\section*{Acknowledgments}
This research was partially supported by NSF awards 2111688,
2450377, 2450378.

ChatGPT (OpenAI) was used for partial assistance with developing code for the numerical solvers, which the authors thoroughly reviewed and verified. All theoretical and numerical results are the authors' own.

\appendix
\renewcommand\theequation{\thesection\arabic{equation}}
\setcounter{equation}{0}

\section{Iterative Approaches for Solving the Necessary Optimality Systems}\label{sec:sinkhorn-all}

\subsection{Cucker--Smale with Phase-Space Endpoint Marginals}\label{sec:sinkhorn-cs-full}
To solve the system in \eqref{eq:Schr-cs-full}, we propose an iterative approach that relies on delayed evaluations of the nonlinear terms. As a first step, we solve the non-interacting second-order Schr\"odinger bridge problem (whose solution is obtained by setting $a = 0$ in \eqref{eq:Schr-cs-full}) using a Sinkhorn iteration similar to \eqref{eq:sinkhorn}. Here, rather than the Fokker--Planck equation \eqref{eq:forward} and its adjoint \eqref{eq:backward}, the Sinkhorn iteration involves the Klein--Kramers equation and its adjoint, namely
\begin{align}
   \partial_t \hat\phi_t +  \big \langle  v , \nabla_x \hat\phi_t  \big \rangle  -\frac{\sigma^2}{2} \Delta_v \hat\phi_t &=0,  \label{eq:forward-kinetic}\\
        \partial_t \phi_t + \big \langle v ,  \nabla_x \phi_t \big \rangle    +\frac{\sigma^2}{2} \Delta_v \phi_t &= 0. \label{eq:backward-kinetic} 
\end{align}
Schematically, the Sinkhorn iteration is as follows.
\begin{align}\label{eq:sinkhorn-kinetic}
\phi_0&(x,v)\qquad \xleftarrow{\eqref{eq:backward-kinetic}} \qquad \! \phi_T(x,v) \nonumber \\
\tfrac{\mu_{\rm initial}(x,v)}{ \phi_0(x,v)}& \bigg\downarrow\quad\qquad\qquad\qquad \qquad\bigg\uparrow\ \tfrac{\mu_{\rm final}(x,v)}{\hat\phi_T(x,v)}\\
\hat \phi&_0(x,v)\qquad \!\!\xrightarrow{\eqref{eq:forward-kinetic}}  \qquad \! \hat \phi_T(x,v).\nonumber
\end{align}
The converged values are used to initialize our iterative scheme by setting
\begin{align*}
\hat \varphi_t^{(0)} = \hat \phi_t,\qquad  \varphi_t^{(0)} =  \phi_t, \qquad  \mu_t^{(0)} = \phi_t\hat \phi_t.
\end{align*}

After obtaining these initializations, we proceed with a scheme consisting of two nested fixed-point iterations.  In the inner iteration, indexed below by $[j]$, Eqs.~\eqref{eq:backward-cs-full} and \eqref{eq:forward-cs-full} are solved with the density $\mu_t$, and hence the nonlocal drift $f_t[\mu_t]$, held fixed. Additionally, this inner iteration also involves a delayed evaluation of the reaction rate $r_t$. 
These delayed evaluations are performed to handle the nonlinear terms, making the numerical solution of the forward-backward equations more 
tractable. 
The outer iteration, indexed below by $(k)$, then updates the density $\mu_t$, which in turn modifies the dependent drift term $f_t[\mu_t]$.

The nested scheme involves computing the values of $\varphi_t, \hat \varphi_t$, and $\mu_t$, all of which are positive functions. We therefore take the \emph{Hilbert metric} denoted by $d_H(\cdot, \cdot)$ below, to define a convergence criterion \cite{georgiou2015positive,bushell1973hilbert}. The Hilbert metric between two positive functions $\varrho, \vartheta$ on the computational phase space is defined as
\begin{align}\label{eq:Hilbert-metric-def}
 \!\! \!\!  d_H(\varrho,\vartheta) := \log \sup_{x,v}\frac{\varrho(x,v)}{\vartheta(x,v)}- \log \inf_{x,v} \frac{\varrho(x,v)}{\vartheta(x,v)}.
\end{align}
Our proposed method is detailed in the following steps. 
\begin{enumerate}
\item \textbf{(Outer fixed-point update)}
At the $(k)$-th outer iteration, the density $\mu_\cdot^{(k)}$ is fixed in the mean-field acceleration, computed as
\begin{align}
     f_\cdot[\mu_\cdot^{(k)}]&(x,v) \nonumber \\
     &= \int a(\Vert x- \tilde{x}\Vert) \mu_\cdot^{(k)}(\tilde{x}, \tilde{v})    \big(\tilde{v} - v \big) \, {\rm d}\tilde{x} \, {\rm d} \tilde{v}.
\end{align}
The pair $(\varphi_\cdot, \hat \varphi_\cdot)$ is initialized by
\begin{align}
\varphi_\cdot^{[0]}=\varphi_\cdot^{(k)}, \qquad  \hat\varphi_\cdot^{[0]}=\hat\varphi_\cdot^{(k)}.
\end{align}

\item \textbf{(Reaction rate)}
Given the function pair \((\varphi_\cdot^{[j]},\hat\varphi_\cdot^{[j]})\), compute the nonlocal reaction rate as
\begin{align}
    r_\cdot^{[j]}&(x,v) = \nonumber \\
    & \int   \hat \varphi_\cdot^{[j]}(\tilde{x},\tilde{v}) a(\Vert x-\tilde{x} \Vert) \big \langle \tilde{v} - v,\nabla_{\tilde{v}} \varphi_\cdot^{[j]}(\tilde{x},\tilde{v}) \big \rangle  \, {\rm d}\tilde{x} \, {\rm d}\tilde{v}.
\end{align}
\item \textbf{(Backward propagation of $\varphi_t$)}
With $r_\cdot^{[j]}$ and $f_\cdot[\mu_\cdot^{(k)}]$ fixed, the next iterate $\varphi_\cdot^{[j+1]}$ is obtained from
\begin{align}\label{eq:backward-alg-cs-full}
   \partial_t \varphi_t^{[j+1]}=  &- \big \langle v , \nabla_x  \varphi_t^{[j+1]} \big \rangle -  \big \langle f_t[\mu_t^{(k)}], \nabla_v \varphi_t^{[j+1]} \big \rangle   \nonumber \\
   &-\frac{\sigma^2}{2} \Delta_v  \varphi_t^{[j+1]} +r_t^{[j]}(x,v)  \varphi_t^{[j+1]},
\end{align}
integrated backward in time with the terminal value
$$\varphi_T^{[j+1]}=\varphi_T^{[j]}.$$
\item \textbf{(Enforcing the initial marginal constraint)}
The initial time constraint is imposed by setting
\begin{align}\label{eq:initial-marginal-alg}
      \hat \varphi_0^{[j+1]} (x,v) =  \frac{\mu_{\rm initial}(x,v)}{\varphi_0^{[j+1]} (x,v)}.
\end{align}
\item \textbf{(Forward propagation of $\hat \varphi_t$)}
With $r_\cdot^{[j]}$ and $f_\cdot[\mu_\cdot^{(k)}]$ fixed, the next iterate $\hat \varphi_\cdot^{[j+1]}$ is obtained from
 \begin{align}\label{eq:forward-alg-cs-full}
          \partial_t \hat\varphi_t^{[j+1]}  = & - \big \langle v, \nabla_x \hat\varphi_t^{[j+1]} \big \rangle  - \nabla_v \cdot(\hat\varphi_t^{[j+1]} f_t[\mu_t^{(k)}]) \nonumber \\
          &+\frac{\sigma^2}{2} \Delta_v \hat\varphi_t^{[j+1]} -r_t^{[j]}(x,v)  \hat\varphi_t^{[j+1]},
    \end{align}
integrated forward in time with the initial value given in \eqref{eq:initial-marginal-alg}.
\item \textbf{(Enforcing the final marginal constraint)}
The final marginal constraint is imposed by setting
\begin{align}
    \varphi_T^{[j+1]}(x,v) =    \frac{\mu_{\rm final}(x,v)}{\hat\varphi_T^{[j+1]}(x,v)} .
\end{align}

\item \textbf{(Inner convergence)}
The inner fixed-point iteration is continued until
\begin{align*}
   \sup_{t\in [0,T]}d_H(\varphi_t^{[j+1]},\varphi_t^{[j]}) & <\epsilon \\
    \sup_{t\in [0,T]}d_H(\hat \varphi_t^{[j+1]},\hat \varphi_t^{[j]}) & < \epsilon,
\end{align*}
where $\epsilon>0$ denotes a prescribed tolerance. 
\item \textbf{(Update of the pair $(\varphi_t, \hat \varphi_t)$)}
Once the inner iteration has converged, the output pair is assigned to the next outer iterates, i.e.,
\begin{align}
    \varphi_\cdot^{(k+1)}=\varphi_\cdot^{[j+1]},    \qquad    \hat\varphi_\cdot^{(k+1)}=\hat\varphi_\cdot^{[j+1]}.\end{align}
\item \textbf{(Relaxed density update)} 
The density $\mu_\cdot$ is updated according to
\begin{align} \label{eq:iteration-nl}
    \mu_\cdot^{(k+1)}  = \theta \varphi_\cdot^{(k+1)} \hat \varphi_\cdot^{(k+1)} + (1-\theta) \mu_\cdot^{(k)}, 
     \end{align}
     for fixed $\theta \in (0,1]$.
\item \textbf{(Outer convergence)}
The outer fixed-point iteration is continued until
\begin{align}
 \sup_{t\in [0,T]}   d_H \big( \mu_t^{(k+1)},\mu_t^{(k)} \big) < \epsilon.
\end{align}
\end{enumerate}

Observe that when the drift and the reaction rate are fixed, \eqref{eq:backward-alg-cs-full} and \eqref{eq:forward-alg-cs-full} become linear PDEs that are numerically more tractable to solve than their nonlinear counterparts, \eqref{eq:backward-cs-full} and \eqref{eq:forward-cs-full}.
Additionally, similar to \eqref{eq:backward-cs-full} and \eqref{eq:forward-cs-full}, the differential operators for \eqref{eq:backward-alg-cs-full} and \eqref{eq:forward-alg-cs-full} are adjoint to each other when the drift and reaction rates are fixed. This adjointness property simplifies the numerical integration of the PDEs; see Appendix~\ref{sec:numerics}. 

The outer iteration consists of a relaxed fixed-point update of $\mu_\cdot$ in \eqref{eq:iteration-nl}, using a fixed damping parameter $\theta \in (0,1]$. For $\theta <1$, this relaxation (known as damping) tempers overly aggressive updates, which can improve the numerical stability \cite{krasnosel1955two,mann1953mean}, whenever necessary.

\subsection{Cucker--Smale with Spatial Endpoint Marginals}\label{sec:sinkhorn-cs-partial}
We adopt the approach detailed in Appendix~\ref{sec:sinkhorn-cs-full} with the following modifications applied to steps $3$ through $6$ of the inner iteration.
\begin{enumerate}[start=3]
\item \textbf{(Backward propagation of $\varphi_t$)}
With $r_\cdot^{[j]}$ and $f_\cdot[\mu_\cdot^{(k)}]$ fixed, the next iterate $\varphi_\cdot^{[j+1]}$ is obtained from
\begin{align}\label{eq:backward-alg-cs-partial}
   \partial_t \varphi_t^{[j+1]}=  &- \big \langle v , \nabla_x  \varphi_t^{[j+1]} \big \rangle -  \big \langle f_t[\mu_t^{(k)}], \nabla_v \varphi_t^{[j+1]} \big \rangle   \nonumber \\
   &-\frac{\sigma^2}{2} \Delta_v  \varphi_t^{[j+1]} +r_t^{[j]}(x,v)  \varphi_t^{[j+1]},
\end{align}
integrated backward in time with the terminal value $$\varphi_T^{[j+1]}(x,v)=\eta_T^{[j]}(x)$$
for all $v$ in the computational domain.
\item \textbf{(Enforcing the initial marginal constraint)}
The initial marginal constraint is imposed by setting
\begin{align}\label{eq:initial-marginal-alg-cs-partial}
      \hat \eta_0^{[j+1]} (x) =   \frac{\rho_{\rm initial}(x)} {\int \nu_0(x,v)  \varphi_0^{[j+1]}(x,v) \,{\rm d} v}.
\end{align}
\item \textbf{(Forward propagation of $\hat \varphi_t$)}
With $r_\cdot^{[j]}$ and $f_\cdot[\mu_\cdot^{(k)}]$ fixed, the next iterate $\hat \varphi_\cdot^{[j+1]}$ is obtained from
 \begin{align}\label{eq:forward-alg-cs-partial}
          \partial_t \hat\varphi_t^{[j+1]}  = & - \big \langle v, \nabla_x \hat\varphi_t^{[j+1]} \big \rangle  -  \nabla_v \cdot(\hat\varphi_t^{[j+1]} f_t[\mu_t^{(k)}]) \nonumber \\
          &+\frac{\sigma^2}{2} \Delta_v \hat\varphi_t^{[j+1]} -r_t^{[j]}(x,v)  \hat\varphi_t^{[j+1]},
    \end{align}
integrated forward in time with the initial value given by
\begin{align}
    \hat \varphi_0^{[j+1]}(x,v) = \hat \eta_0^{[j+1]}(x)\, \nu_0(x,v).
\end{align}
\item \textbf{(Enforcing the final marginal constraint)}
The terminal marginal constraint is imposed by setting
\begin{align}
    \eta_T^{[j+1]} (x) &=  \frac{\rho_{\rm final}(x)}{\int \hat\varphi_{T}^{[j+1]}(x,v)\, {\rm d} v}.
\end{align}
\end{enumerate}

The iterative scheme can be initialized using a non-interacting Schr\"odinger bridge, corresponding to $a = 0$. In the present setting, this bridge satisfies the spatial endpoint-marginal constraints and is characterized by the pair $(\phi_\cdot, \hat \phi_\cdot)$, which solve
\begin{subequations}
       \begin{align}
 \partial_t \hat\phi_t + \big \langle  v , \nabla_x \hat\phi_t  \big \rangle  -\frac{\sigma^2}{2} \Delta_v \hat\phi_t &=0,  \label{eq:forward-kinetic-2}\\
 \partial_t \phi_t + \big \langle v ,  \nabla_x \phi_t \big \rangle    +\frac{\sigma^2}{2} \Delta_v \phi_t &= 0, \label{eq:backward-kinetic-2} 
\end{align}
subject to the same conditions in \eqref{eq:BCs-cs-partial-all-1}-\eqref{eq:BCs-cs-partial-all-2}.
\end{subequations}
Numerical solutions can be obtained via a Sinkhorn iteration:
\begin{align*}
\phi_0&(x,v)\quad \xleftarrow{\eqref{eq:backward-kinetic-2}} \quad \! \phi_T(x,v) \nonumber \\
\tfrac{\rho_{\rm initial}(x)\nu_0(x,v)}{\int \nu_0(x,v)\phi_0(x,v) \, {\rm d}v}& \bigg\downarrow\quad\qquad\qquad\quad \qquad\bigg\uparrow\ \tfrac{\rho_{\rm final}(x)}{\int \hat\phi_T(x,v) \, {\rm d} v}\\
\hat \phi&_0(x,v)\quad \!\!\xrightarrow{\eqref{eq:forward-kinetic-2}}  \quad \! \hat \phi_T(x,v).\nonumber
\end{align*}

\subsection{Morse Interaction with Phase-Space or Spatial Endpoint Marginals}\label{sec:sinkhorn-morse}
The system \eqref{eq:Schr-m-full} can be solved by an approach analogous to that described in Appendix~\ref{sec:sinkhorn-cs-full}. Specifically, Steps $1,2,3,$ and $5$ are modified by replacing the Cucker--Smale drift $f_\cdot[\mu_\cdot^{(k)}]$ and the reaction rate $r_\cdot^{[j]}$ with 
\begin{align}
     g_\cdot[\rho_\cdot^{(k)}](x) = -\int_{\mathbb R^d} \nabla_x W(\Vert x - \tilde{x} \Vert) \rho_\cdot^{(k)}(\tilde{x})\, {\rm d} \tilde{x}\\
    \text{where } \quad  \rho_\cdot^{(k)}(\cdot) = \int \mu_\cdot^{(k)}(\cdot,v) \, {\rm d} v,
\end{align} 
and 
\begin{align}
    &s_\cdot^{[j]}(x) = \nonumber \\
    &-\int \!  \Big \langle\nabla_x W(\Vert x -\tilde{x} \Vert),   \nabla_{\tilde v} \varphi_\cdot^{[j]}(\tilde{x},\tilde{v}) \Big \rangle \hat \varphi_\cdot^{[j]}(\tilde{x},\tilde{v}) \, {\rm d}\tilde{x}\: {\rm d}\tilde{v}. 
\end{align} 
As noted in Section~\ref{sec:morse-}, since $g_t[\rho_t]$ is independent of $v$, the advection term $\nabla_v \cdot (g_t[\rho_t] \hat\varphi_t)$ in the forward equation simplifies, relative to the Cucker--Smale case, to $\langle \nabla_v \hat\varphi_t , g_t[\rho_t] \rangle$. 
The rest of the steps remain the same. 

Likewise, to solve the optimality system of \eqref{problem:morse-partial}, we apply a similar procedure to that in Appendix~\ref{sec:sinkhorn-cs-partial} but adapted for $g_\cdot[\rho_\cdot^{(k)}]$ and $s_{\cdot}^{[j]}$.

\section{Numerical Solutions of the Vlasov--Fokker--Planck Equations} \label{sec:numerics}
This appendix discusses how to integrate the backward and forward equations in the iterative schemes from Appendix~\ref{sec:sinkhorn-all}. We exploit the adjointness to simplify the computation. The exposition here focuses on Eqs.~\eqref{eq:backward-alg-cs-full},\eqref{eq:forward-alg-cs-full} in Appendix~\ref{sec:sinkhorn-cs-full}; however, the approach applies analogously to the adjoint equations involved in Appendices~\ref {sec:sinkhorn-cs-partial} and \ref{sec:sinkhorn-morse}.

\subsection{Adjointness in the Optimality Systems}\label{sec:adjoint}
Let $\mathcal A_t$ denote the advection operator
\begin{align*}
  \mathcal A_t  h :=  -  \langle v , \nabla_x h \rangle - \langle f_t[\mu_t] ,\nabla_v h \rangle.
\end{align*}
Then we can write \eqref{eq:backward-alg-cs-full},\eqref{eq:forward-alg-cs-full} as
\begin{align}
\partial_t  \varphi_t&= \underbrace{\Big(\mathcal A_t  - \frac{\sigma^2}{2} \Delta_v +r_t \Big)}_{\displaystyle =: \mathcal B_t} \varphi_t, \label{eq:diff-phij} \\
   \partial_t \hat \varphi_t &= \underbrace{\Big(\mathcal A_t  - (\nabla_v \cdot f_t[\mu_t]) +\frac{\sigma^2}{2} \Delta_v -r_t \Big)}_{\displaystyle =: \mathcal F_t} \hat \varphi_t, \label{eq:diff-hatphij}
\end{align}
where we have dropped the iteration indices for brevity. In \eqref{eq:diff-hatphij}, the quantities  $(\nabla_v \cdot f_t[\mu_t])$ and $r_t$ are viewed as operators that act via multiplication. Here, $\mathcal B_t$ and $\mathcal F_t$ represent the backward and forward differential operators, respectively.  

We define the adjoint of $\mathcal B_t$, denoted as $\mathcal B_t^\dagger$, via
\begin{align}
    \big \langle \mathcal B_t h_1, h_2 \big \rangle_{L^2} = \big \langle h_1 , \mathcal B_t^\dagger h_2 \big \rangle_{L^2},
\end{align}
where $\langle \cdot, \cdot \rangle_{L^2}$ denotes the $L^2$ inner product.
From \eqref{eq:diff-phij}, we have
\begin{align}
    \mathcal B_t^\dagger = \mathcal A_t^\dagger - \frac{\sigma^2}{2} \Delta_v^\dagger +r_t^\dagger.
\end{align}
We can calculate all such adjoint operators using integration by parts:
\begin{align*}
 \big \langle \mathcal A_t h_1 ,& h_2\big \rangle_{L^2}  \\
 &=  \int  -\big \langle v , \nabla_x h_1 \big \rangle  h_2-\big \langle f_t[\mu_t] , \nabla_v h_1 \big \rangle  h_2 \, {\rm d} x \, {\rm d} v\\
    &=  \int   h_1\big \langle v , \nabla_x h_2  \big \rangle  + h_1\nabla_v \cdot (f_t[\mu_t] h_2)  \, {\rm d} x \, {\rm d} v \\
    &=  \big \langle h_1 ,  \mathcal A_t^\dagger h_2\big \rangle_{L^2},
    \end{align*}
    with $$\mathcal A_t^\dagger h = v \cdot \nabla_x h+ \nabla_v \cdot (f_t[\mu_t]h ).$$
    Similarly, one can verify that $\Delta_v^\dagger = \Delta_v$ and $r_t^\dagger = r_t$. It follows that
\begin{align}
    \mathcal B_t^\dagger &= -\mathcal F_t, \label{eq:adj-A-B} \\
    \big \langle \mathcal B_t \varphi_t, \hat \varphi_t \big \rangle_{L^2} & = -  \big \langle  \varphi_t, \mathcal F_t \hat \varphi_t\big \rangle_{L^2}.\label{eq:adj-A-B-with-phis}
\end{align}

Consider a time interval $[t,t+\delta t]$ where $\delta t$ is sufficiently small. Let $\Pi_{t,t+\delta t}$ and $\hat \Pi_{t,t+\delta t}$ denote the backward propagator and forward propagator for $\varphi_{t+\delta t}$ and $\hat \varphi_t$, respectively, over $[t, t+\delta t]$. That is,
\begin{align}
    \varphi_{t} &= \Pi_{t,t+\delta t} \, \varphi_{t+\delta t} =e^{-\delta t \mathcal B_t} \varphi_{t+\delta t}, \label{eq:backward-propagator} \\
    \hat \varphi_{t+\delta t} &=  \hat  \Pi_{t,t+\delta t} \, \hat \varphi_{t} = e^{\delta t \mathcal F_t} \, \hat \varphi_{t},  \label{eq:forward-propagator}
\end{align}
where the operators $\mathcal B_t$ and $\mathcal F_t$ can be considered fixed over the indicated small enough interval. 
By \eqref{eq:adj-A-B}, we have
\begin{align}
\hat \Pi_{t,t+\delta t}     
&= e^{-\delta t \mathcal B_t^\dagger}  \nonumber \\
&= I-\delta t \mathcal B_t^\dagger + \frac{\delta t^2 (\mathcal B_t^\dagger)^2}{2} - \ldots \nonumber \\
& =  (e^{-\delta t \mathcal B_t})^\dagger \label{eq:dagger-out} \\
& =  \Pi_{t,t+\delta t}^\dagger. \label{eq:discrete-operator-dagger}
\end{align}
Eq.~\eqref{eq:dagger-out} follows from the fact that $(A_1 A_2)^\dagger = A_2^\dagger A_1^\dagger$ and Eq.~\eqref{eq:discrete-operator-dagger} attests that the propagators are adjoint.

In a discretized domain, $\varphi_t,\hat \varphi_t$ are vectors, and the propagators are matrices. Moreover, the $L^2$ inner product can be expressed in the discrete setting as
\begin{align*}
     \langle \varphi_t,  \hat \varphi_t \rangle_{{\rm dis}} = \varphi_t^\top \hat \varphi_t  \, (\delta x)^d \, (\delta v)^d,
\end{align*}
assuming a uniform grid. Therefore, by \eqref{eq:discrete-operator-dagger}, we have
\begin{align}
 \langle \Pi_{t, t+\delta t} \varphi_{t+\delta t},  \hat \varphi_t \rangle_{{\rm dis}} =
 \langle  \varphi_{t+\delta t}, \hat \Pi_{t, t+\delta t} \hat \varphi_t \rangle_{{\rm dis}}.\label{eq:approx-1}
\end{align}
Eq.~\eqref{eq:approx-1} implies
\begin{align*}
    \varphi_{t+\delta t}^\top (\Pi_{t,t+\delta t}^\top - \hat\Pi_{t,t+\delta t}) \hat  \varphi_t = 0,
\end{align*}
for any such $\hat  \varphi_t,\varphi_{t+\delta t}^\top$. Hence,
\begin{align}
    \Pi_{t,t+\delta t}^\top = \hat\Pi_{t,t+\delta t}.
    \end{align}
That is, the matrix $\Pi_{t,t+\delta t}$ propagates $\varphi_{t+\delta t}$ backward in time, and its transpose propagates $\hat \varphi_{t}$ forward in time. As such, it suffices to compute one of these propagators and obtain the other via taking the transpose. We next expand on the computation of the discrete forward propagator $\hat \Pi_{t , t+\delta t}$.

\subsection{Strang Splitting}

We use Strang splitting \cite{strang1968construction} for approximating the discrete forward propagator that corresponds to a Vlasov--Fokker--Planck equation, as in \eqref{eq:diff-hatphij}, described in \eqref{eq:forward-propagator}. This method has been used in related works (see, e.g., \cite{carrillo2019monte,joglekar2020vlapy}).

Looking more closely at \eqref{eq:diff-hatphij}, we notice the following structure.
\begin{widetext}
\begin{align}
\partial_t \hat \varphi_t
= \underbrace{- \langle v , \nabla_x \hat\varphi_t \rangle- \langle   f_t[\mu_t], \nabla_v \hat\varphi_t\rangle \vphantom{\frac{\sigma^2}{2}\big[\nabla_v \cdot f_t[\mu_t]+r_t\big]}
}_{\text{Advection}} \underbrace{+ \frac{\sigma^2}{2}\Delta_v \hat\varphi_t \vphantom{\frac{\sigma^2}{2}\big[\nabla_v \cdot f_t[\mu_t]+r_t\big]}}_{\text{Diffusion}}
\underbrace{ -\big(\nabla_v \cdot f_t[\mu_t]+r_t\big)\hat\varphi_t\vphantom{\frac{\sigma^2}{2}\big[\nabla_v \cdot f_t[\mu_t]+r_t\big]}
}_{\text{Reaction}}.
\end{align}
\end{widetext}
For a discretized domain, we rewrite this equation as
\begin{align}
    \partial_t \hat \varphi_t = (\mathcal F_1  + \mathcal F_2+ \mathcal F_3) \hat \varphi_t,
\end{align}
where $\mathcal F_1$, $\mathcal F_2$, and $\mathcal F_3$ represent discrete approximations of advection, diffusion, and reaction operators, respectively. Strang Splitting considers the second-order approximation of solutions of the previous equation via
\begin{align}\label{eq:forward-strang}
      \hat \varphi_{t+\delta t} = (e^{\frac{\delta t}{2}\mathcal F_3}) (e^{\frac{\delta t}{2}\mathcal F_2}) (e^{\delta t \mathcal F_1})  (e^{\frac{\delta t}{2}\mathcal F_2}) (e^{\frac{\delta t}{2}\mathcal F_3}) \hat \varphi_t + O(\delta t^3).
\end{align}
Both $\mathcal F_1$ and $\mathcal F_3$ are time-dependent; 
however, over a sufficiently small interval $[t, t+\delta t]$, they can be held fixed by evaluating them at the midpoint $t+ (\delta t/2)$ with the midpoint values approximated by averaging their values at the interval endpoints. The Strang splitting procedure is as follows. Starting from $\hat \varphi_t$, 
\begin{enumerate}
    \item compute $\zeta_1 = e^{\frac{\delta t}{2}\mathcal F_3} \hat \varphi_t$
    \item compute $\zeta_2 = e^{\frac{\delta t}{2}\mathcal F_2} \zeta_1$
    \item compute $\zeta_3 = e^{\delta t\mathcal F_1} \zeta_2$\\[2.5pt]
    \item compute $\zeta_4 =  e^{\frac{\delta t}{2}\mathcal F_2} \zeta_3$
    \item  compute $\hat \varphi_{t+\delta t} =  e^{\frac{\delta t}{2}\mathcal F_3} \zeta_4$.
\end{enumerate}
The advection step is computed using a semi-Lagrangian method \cite{cheng1976integration} so as to obtain characteristic trajectories. This can be done using, e.g., an explicit Runge-Kutta solver or a Verlet symplectic scheme \cite[Chapter 8]{sarkka2019applied}. The Verlet scheme is suitable for Hamiltonian vector fields, such as $\begin{bmatrix}v, g_t[\rho_t](x)\end{bmatrix}^\top$ in Morse-potential-driven kinetics. Eq.~\eqref{eq:forward-strang} gives an approximation of the forward propagator $\hat \Pi_{t, t+\delta t}$. As explained in Appendix~\ref{sec:adjoint}, the backward propagator can be approximated by taking the transpose.

\bibliography{Paper2-Uploads/references}

\end{document}